\RequirePackage{snapshot}

\documentclass{siamltex}

\usepackage[text={5.125in,8.25in},centering]{geometry}
\usepackage{microtype}
\usepackage{xspace}      %
\usepackage{amsmath}
\usepackage{amssymb}
\usepackage{enumitem}
\usepackage{overpic}
\usepackage{hyperref}
\usepackage{graphicx}
\usepackage{cite}        %
\usepackage{braket}
\usepackage{subdepth}
\usepackage{color}

\usepackage{tikzexternal}
\newtheorem{defn}[theorem]{Definition}
\def\T{^T\!}

\newcommand{\norm}[1]{\left\|#1\right\|}
\newcommand{\tnorm}[1]{\norm{#1}_2}

\newcommand{\ip}[2]{\left\langle #1,\, #2\right\rangle}
\newcommand{\map}[3]{#1: #2\rightarrow #3}

\let\Newset\Set
\renewcommand{\set}[2]{\Newset{#1|#2}}
\newcommand{\dom}[1]{\mathop{\mathrm{dom}} #1}
\newcommand{\cl}[1]{\mathrm{cl}\left(#1\right)}
\newcommand{\argmin}{\mathop{\mathrm{arg\,min}}}
\newcommand{\argmax}{\mathop{\mathrm{arg\,max}}}

\newcommand{\uball}{I\!\!B}

\newcommand{\half}{{\textstyle \frac{1}{2}}}

\newcommand{\support}[2]{\delta^*\!\left(#1 \left| \vphantom{#1}\, #2 \right.\right)}
\newcommand{\indicator}[2]{\delta\left(#1 \mid #2 \right)}
\newcommand{\ncone}[2]{N\left(#1 \left|\vphantom{#1}\, #2\right.\right)}

\newcommand{\lev}[2]{\mathrm{lev}_{#1}(#2)}
\newcommand{\sign}{\mathrm{sign}}
\newcommand{\gauge}[2]{\gamma\left(#1 \mid #2 \right)}
\newcommand{\cone}[1]{\mathrm{cone}\left(#1\right)}

\newcommand{\intr}[1]{\mathrm{int}\left(#1\right)}
\newcommand{\epi}[1]{\mathrm{epi}\ #1}

\newcommand{\Span}[1]{\mathrm{span}\left(#1\right)}
\newcommand{\ri}[1]{\mathrm{ri}\left(#1\right)}
\newcommand{\infc}{{\triangledown}}

\newcommand{\barrier}[1]{\mathrm{bar}\left(#1\right)}
\newcommand{\sd}{\partial}
\newcommand{\hzn}[1]{\mathrm{hzn}\left(#1\right)}
\newcommand{\barr}[1]{\mathrm{bar}\left(#1\right)}

\newcommand{\lam}{\lambda}

\newcommand{\eps}{\epsilon}

\newcommand{\alf}{\alpha}
\newcommand{\del}{\delta}

\newcommand{\bR}{\mathbb{R}}

\newcommand{\R}{\bR}
\newcommand{\Rn}{\R^n}
\newcommand{\Rnn}{\R^{n\times n}}
\newcommand{\Rm}{\R^m}
\newcommand{\Rmn}{\R^{m\times n}}

\newcommand{\bB}{\mathbb{B}}
\newcommand{\eR}{{\overline\bR}}
\newcommand{\cP}{\mathcal{P}}

\newcommand{\cD}{{\mathcal{D}}}

\newcommand{\cS}{\mathcal{S}}

\newcommand{\bmu}{{\overline \mu}}

\newcommand{\bv}{\overline v}
\newcommand{\bu}{\overline u}

\newcommand{\br}{\overline r}
\newcommand{\bs}{\overline s}

\newcommand{\bw}{\overline w}

\newcommand{\bx}{\overline x}

\newcommand{\hx}{{\hat x}}

\newcommand{\hb}{{\hat b}}
\newcommand{\hu}{{\hat u}}

\newcommand{\hr}{{\hat r}}
\newcommand{\hA}{{\hat A}}

\newcommand{\hrho}{{\hat \rho}}
\newcommand{\hphi}{{\hat \phi}}
\newcommand{\tx}{{\tilde x}}
\newcommand{\ty}{{\tilde y}}
\newcommand{\tb}{{\tilde b}}

\newcommand{\tw}{{\tilde w}}
\newcommand{\tA}{{\tilde A}}
\newcommand{\tB}{{\tilde B}}
\newcommand{\tU}{{\tilde U}}
\newcommand{\trho}{{\tilde \rho}}
\newcommand{\tphi}{{\tilde \phi}}

\newcommand{\ssL}{{\mbox{\tiny L}}}
\newcommand{\ssR}{{\mbox{\tiny R}}}

\newcommand{\hv}{\hat v}

\newcommand{\RQP}{QS\xspace}
\newcommand{\Zal}{Z$\breve{\mathrm{a}}$linescu\xspace}

\newcommand{\B}[1]{{\bf #1}}
\newcommand{\Ran}[1]{\mathrm{Ran}\left(#1\right)}
\newcommand{\Nul}[1]{\mathrm{Nul}\left(#1\right)}
\renewcommand{\text}[1]{\hbox{\quad#1\quad}}
\newcommand{\Pbt}{\ensuremath{\cP(b,\tau)}\xspace}
\newcommand{\Pbl}{\ensuremath{\cP_\ssL(b,\lambda)}\xspace}
\newcommand{\Pbs}{\ensuremath{\cP_\ssR(b,\sigma)}\xspace}

\def\minim{\mathop{\hbox{\rm minimize}}}
\def\maxim{\mathop{\hbox{\rm maximize}}}
\def\minimize#1{\displaystyle\minim_{#1}}

\def\st{\mathop{\hbox{\rm subject to}}}

\definecolor{softgray}{rgb}{0.92,0.92,0.95}
\definecolor{softblue}{rgb}{0.90,0.92,1.00}
\definecolor{lightgray}{rgb}{0.12,0.12,0.55}

\definecolor{theframe}{gray}{0.75}
\definecolor{theblue} {rgb}{0.02,0.04,0.48}
\definecolor{thegrey} {gray}{0.5}
\definecolor{theshade}{gray}{0.98}
\definecolor{thered}  {rgb}{0.00,0.00,0.00}
\definecolor{thegreen}{rgb}{0.3,0.3,0.3}
\newcommand{\spgl}{SPGL1\xspace}

\usepackage{color}
\definecolor{softblue}{rgb}{0.90,0.92,1.00}
\usepackage{mdframed}
\mdfdefinestyle{mystyle}{%
    outerlinewidth  = 2             ,%
    leftmargin      = 40            ,%
    rightmargin     = 40            ,%
    backgroundcolor = softblue      ,%
    outerlinecolor  = blue          ,%
    innertopmargin  = .5\topskip      ,%
    splittopskip    = 1.1\topskip   ,%
    splitbottomskip = 0.5\topskip,
    ntheorem        = true          ,%
    skipabove       = \baselineskip ,%
    skipbelow       = \baselineskip]
}
\surroundwithmdframed[style=mystyle]{theorem}
\surroundwithmdframed[style=mystyle]{lemma}
\surroundwithmdframed[style=mystyle]{corollary}

  \title{Variational Properties of Value Functions}
  \author{%
    Aleksandr Y. Aravkin\thanks{IBM T.J. Watson Research Center;
      (saravkin@us.ibm.com)} The work of this author was
    partially supported by an NSERC CRD grant through the University of British Columbia.
    \and James V. Burke\thanks{Department of Mathematics, University
      of Washington, Seattle, WA 98195 (burke@math.washington.edu).}
    \and Michael P. Friedlander\thanks{Department of Computer
      Science, University of British Columbia, Vancouver B.C., V6T 1Z4,
      Canada (mpf@cs.ubc.ca). The work of this author was
      partially supported by NSERC Discovery Grant 312104.%
    \hfill November 15, 2012; revised May 9, 2013; revised May 23, 2013}}
\begin{document}

\maketitle

\begin{abstract}
Regularization plays a key role in a variety of optimization formulations of inverse problems. A recurring theme in regularization approaches is the selection of regularization parameters, and their effect on the solution and on the optimal value of the optimization problem.  The sensitivity of the value function to the regularization parameter can be linked directly to the Lagrange multipliers.  This paper characterizes the variational properties of the value functions for a broad class of convex formulations, which are not all covered by standard Lagrange multiplier theory. An inverse function theorem is given that links the value functions of different regularization formulations (not necessarily convex).  These results have implications for the selection of regularization parameters, and the development of specialized algorithms. Numerical examples illustrate the theoretical results.
\end{abstract}


\section{Introduction} \label{sec:introduction}

It is well known that there is a close connection between the
sensitivity of the optimal value of a parametric optimization problem
and its Lagrange multipliers. Consider the family of feasible convex
optimization problems
\begin{equation*}
  \tag*{\Pbt}
  \minimize{r, x} \quad \rho(r) \text{subject to} Ax + r = b, \quad \phi(x) \le \tau,
\end{equation*}
where $b\in\Rm$, $A\in\Rmn$, and the functions
$\map{\phi}{\Rn}{\eR:=(-\infty, \infty]}$ and $\map{\rho}{\Rm}{\eR}$
are closed, proper and convex, and continuous relative to their
domains. %
The value function
\begin{equation*}
  v(b,r) := \inf_{r,x}\set{ \rho(r)}{Ax+r=b,\ \phi(x) \leq \tau}
\end{equation*}
gives the optimal objective value of problem $\Pbt$ for fixed
parameters $b$ and $\tau$. If $\Pbt$ is a feasible {\it ordinary
  convex program}~\cite[Section 28]{RTR}, then under standard
hypotheses the subdifferential of $v$ is
the set of pairs $(u,\mu)$,
where $u\in\Rm$ and $\mu\in\R$ are the Lagrange multipliers of $\Pbt$
corresponding to the equality and inequality constraints, respectively. This
connection is extensively explored in Rockafeller's 1993 survey
paper~\cite{Rockafellar:1993}.

If we allow $\phi$ to take on infinite values on the domain of the 
objective---which can occur, for example, if $\phi$ is an arbitrary gauge---%
then \Pbt is no longer an ordinary convex program, and so the standard
Lagrange multiplier theory does not apply.  Multiplier theories that do
apply to more general contexts can be found in
\cite{EkeTem:1976,ComLagThi:1994,Zal:2002,BauCom:2011}.
Remarkably, even in this general setting, it is possible to obtain 
explicit formulas of the subdifferential of the value function $v$ 
useful in many applications. %

\subsection{Examples}
We give two simple examples that illustrate the need for the extended
Lagrange multiplier theory. Both are of the form
\begin{equation}
  \label{LittleLS}
  \minimize{x} \quad
  \half \tnorm{Ax-b}^2
  \text{subject to}
  \gauge{x}{U}\le 1,
\end{equation}
where
\[
\gauge{x}{U}:=\inf\set{\lambda\ge0 }{ x\in \lambda U}
\] 
is the {\it gauge function} for the closed nonempty convex set
$U\subset\Rn$, which contains 0. Let $A = I$ and $b = (0, -1)^T$.
Then the solution to \eqref{LittleLS} is just the 2-norm projection onto
the set $\set{x}{\gauge{x}{U}\le 1}=U$.

For our first example, we consider the set \[
U=\set{x\in\bR^2}{\half x_1^2\le x_2},%
\]
defined in~\cite[Section 10]{RTR}.
The gauge for this set is an example of a closed, proper and
convex function that is not locally bounded and therefore not continuous at a point in its
effective domain.  It is straightforward to show that
\[
\gauge{x}{U}=
\begin{cases}
\frac{x_1^2}{2x_2},&x_2>0\\
0,&x_1=0=x_2\\
+\infty,&\mbox{otherwise.}\end{cases}
\]
The constraint region for~\eqref{LittleLS} is the set $U$ and
the unique global solution is the point $x=0$. However, 
since $0=\gauge{0}{U}<1$, the classical Lagrange multiplier theory fails: 
the solution is on the boundary of the feasible region, 
and yet no classical Lagrange multiplier exists. 
The problem is that the constraint is active at the solution, but not active 
in the functional sense, i.e., $\gauge{0}{U} < 1$.  
In contrast, the extended multiplier theory of \cite[Theorem 2.9.3]{Zal:2002} %
succeeds with the multiplier choice of $0$.

For the second example, take $U=\bB_2 \cap K$, where $\bB_2$ is the
unit ball associated with the Euclidean norm on $\R^2$.
Then $\gauge{x}{\bB_2\cap
  K}=\tnorm{x}+\indicator{x}{K}$, and the constraint region
for~\eqref{LittleLS} is the set $\bB_2\cap K$. 
Set
$K=\set{(x_1,x_2)}{x_2\ge 0}$. Again, the origin is
the unique global solution to this optimization problem, and no
classical Lagrange multiplier for this problem exists.  

In both of these examples, the multiplier theory in \cite{Zal:2002} 
can be applied to obtain a Lagrange multiplier theorem.
In Theorem \ref{reduced primal-dual}, we extend this theory and provide a 
characterization of these Lagrange multipliers that is useful in computation.

\subsection{Formulations} \label{sec:formulations}

Appropriate definitions of the functions $\rho$ and $\phi$ can be used
to represent a range of practical problems. Choosing $\rho$ to be the
2-norm and $\phi$ to be any norm yields the canonical regularized
least-squares problem
\begin{equation} \label{eq:3}
  \minim_{x,\, r} \quad \|r\|_2
  \quad\st\quad
  \ Ax+r=b,\ \|x\| \le \tau,
\end{equation}
which optimizes the misfit between the data $b$ and the forward model
$Ax$, subject to keeping $x$ appropriately bounded in some norm. The
2-norm constraint on $x$ yields a Tikhonov regularization, popular in
many inversion applications. A 1-norm constraint on $x$ yields the
Lasso problem~\cite{Tib96}, often used in sparse recovery and
model-selection applications.  Interestingly, when the optimal
residual $r$ is nonzero, the value function for this family of
problems is always differentiable in both $b$ and $\tau$
with
\begin{equation*}
  \nabla v(b,\tau) =
  \left(
    \frac{r}{\|r\|_2},\ \frac{\|A^T r\|_*}{\|r\|_2}
  \right),
\end{equation*}
where $\|\cdot\|_*$ is the norm dual to $\|\cdot\|$. This gradient is
derived by van den Berg and
Friedlander~\cite[Theorem~2.2]{BergFriedlander:2008}. The analysis of the sensitivity in
$\tau$ of the value function for the Lasso problem led to the
development of the \spgl solver~\cite{BergFrie:2007b}, currently used
in a variety of sparse inverse problems, with particular success in
large-scale sparse inverse problems~\cite{HFY:2012}. A subsequent
analysis~\cite{BergFriedlander:2011} that allows $\phi(x)$ to be a
gauge paved the way for other applications, such as group-sparsity
promotion~\cite{BergFrie:2010}.

An alternative to $\Pbt$ is the class of penalized formulations
\begin{equation}
  \tag*{\Pbl}
      \minimize{x}\quad \rho(b-Ax) + \lambda\phi(x)
\end{equation}
(the subscript ``$L$'' in the label reminds us that it can be
interpreted as a Lagrangian of the original problem).  The nonnegative
regularization parameter $\lambda$ is used to control the tradeoff
between the data misfit $\rho$ and the regularization term $\phi$. For
example, talking $\rho(r)=\|r\|_2$ and $\phi(x) = \|x\|$ yields a
formulation analogous to~\eqref{eq:3}. This penalized formulation is
commonly used in applications of Bayesian parametric
regression~\cite{MacKay,Roweis99,McKayARD,Tipping2001,Wipf_IEEE_TIT_2011},
inference problems on dynamic linear
systems~\cite{Anderson:1979,Brockett}, feature selection, selective
shrinkage and compressed sensing~\cite{Hastie90,LARS2004,Donoho2006},
robust formulations~\cite{Huber81,Gao2008,Aravkin2011,Farahmand2011},
support-vector regression~\cite{Vapnik98,Hastie01},
classification~\cite{Evgeniou99,Pontil98,Scholkopf00}, and functional
reconstruction~\cite{Aronszajn,Saitoh,Cucker}.

\begin{figure}[t]
  \centering
  \begin{tabular}{@{}c@{}c@{}}
    \includegraphics[width=.48\textwidth]{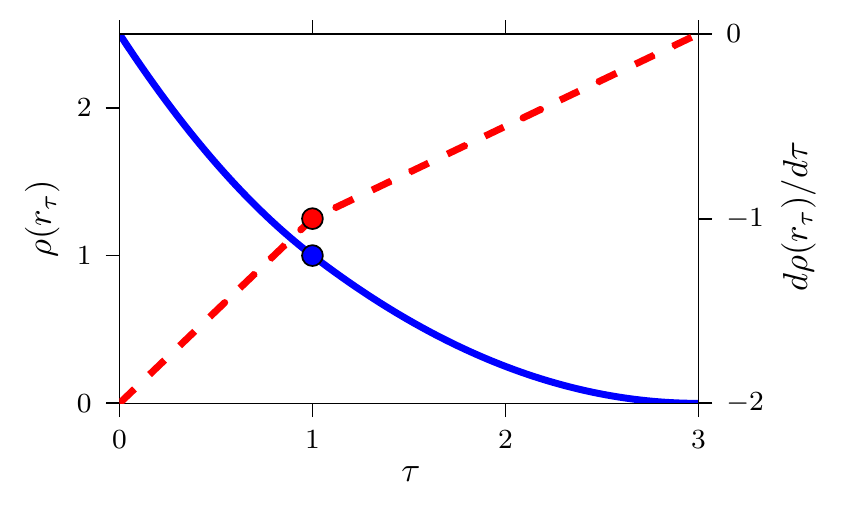}
   &\includegraphics[width=.48\textwidth]{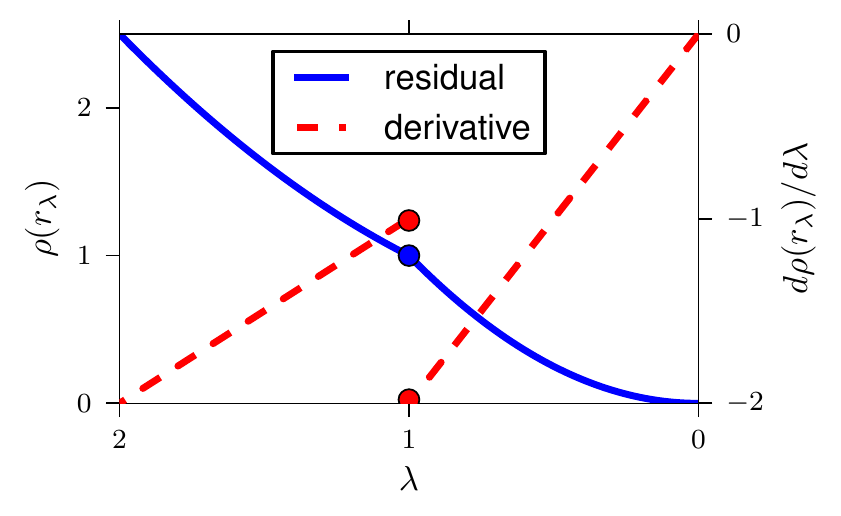}
  \end{tabular}
  \caption{The misfit $\rho(r)$ (solid line) and its derivative
    (dashed line) as a function of the regularization parameter for a
    1-norm regularized example. The left panel shows the constrained
    formulation $\Pbt$, which varies smoothly with $\tau$; the
    right panel shows that the penalized formulation does not vary
    smoothly with $\lambda$ (note the reversed axis).}
  \label{fig:example_tau_vs_rho}
\end{figure}

From an algorithmic point of view, the unconstrained formulation \Pbl
may be preferable. However, the constrained formulation \Pbt has the
distinction that its value function $v(b,\tau)$ is jointly
convex in its parameters; see section~\ref{sec:approach}.  In
contrast, the optimal value of the penalized formulation $\Pbl$ is not
in general a convex function of its parameters. The following simple
example
\[
\rho(r)=\half\|r\|_2^2, \quad \phi(x)=\|x\|_1, \quad
A = \begin{bmatrix}1 & 0 \\ 0 & 1\end{bmatrix}, \quad
b = \begin{bmatrix}2 \\ 1\end{bmatrix},
\]
illustrates this situation.  The optimal values of $\rho$ in the
formulations \Pbt and \Pbl, as functions of $\tau$ and $\lambda$,
respectively, are given by
\[
\rho_\tau =
   \begin{cases}
        \half + \half(\tau-2)^2     & \mbox{for $\tau\in[0,1)$}
      \\\textstyle\frac14(\tau-3)^2 & \mbox{for $\tau\in[1,3)$}
      \\0                           & \mbox{otherwise;}
   \end{cases}
\qquad
\rho_\lambda =
   \begin{cases}
        \lambda^2                   & \mbox{for $\lambda\in[0,1)$}
      \\\half + \half\lambda^2      & \mbox{for $\lambda\in[1,2)$}
      \\5/2                         & \mbox{otherwise.}
   \end{cases}
\]
The optimal values and their derivatives are shown in
Figure~\ref{fig:example_tau_vs_rho}, where it is clear that
$\rho_\tau$ is convex (and in this case also smooth) in $\tau$, but
$\rho_\lambda$ is not convex in $\lambda$.

The admissibility of variational analysis and convexity of the value
function may convince some practitioners to explore formulations of
type $\Pbt$ rather than $\Pbl$. In fact, we give an example
(in section~\ref{Numerics}) of how this variational information can
be used for algorithm design in the context of large-scale inverse
problems.

\subsection{Approach} \label{sec:approach}
 
For many practical inverse problems, the formulation of primary
interest is the residual-constrained formulation
\begin{equation*}
  \tag*{\Pbs}
  \minimize{x}\quad\phi(x)\quad\st\quad\rho(b-Ax)\le\sigma,
\end{equation*}
(the subscript ``R" reminds us that this formulation reverses the
objective and constraint functions from that of \Pbt) in part because
estimates of a tolerance level $\sigma$ on fitting the error
$\rho(b-Ax)$ are more easily available than estimates of a bound on
the penalty parameter on the regularization $\phi$; cf.~\Pbl. However,
the formulation \Pbt can sometimes be easier to solve. The underlying
numerical theme is to develop methods for solving \Pbs
that use a sequence of solutions to the possibly easier problem \Pbt.

In section~\ref{IFT}, we present an inverse function theorem for value
functions that characterizes the relationship between \Pbt and \Pbs,
and applies more generally to nonconvex problems.  
Pairs of problems of this type are classical, though typically 
paired in a max-min fashion. For example, the isoperimetric inequality
and Queen Dido's problem are of this type; the greatest area surrounded by a curve of 
given length is related to the problem of finding the curve of least arc length
surrounding a given area (see \cite{Tap:2013} for a modern survey). 
The Markowitz mean-variance portfolio theory is
also based on such a pairing; minimizing volatility subject to a lower bound on
expected return is related to maximizing expected return subject to an upper
bound on volatility \cite{Mar:1987}.

The application motivating our investigation
is establishing conditions under which it is possible to
implement a root-finding approach for the nonlinear equation
\begin{equation}
  \label{eq:v-tau-sig}
 \mbox{find $\tau$ such that $v(b,\tau) = \sigma,$}
\end{equation}
where \Pbs can be solved via a sequence of approximate solutions of
\Pbt. This generalizes the approach used by van den Berg and
Friedlander~\cite{BergFriedlander:2008,BergFriedlander:2011} for
large-scale sparse optimization applications. The convex case is
especially convenient, because both value functions are decreasing and
convex. When the value function is differentiable, Newton's method is
globally monotonic and locally quadratic. In
section~\ref{VariationalProperties} we establish the variational
properties (including conditions necessary for differentiability) of
\Pbt.

In section~\ref{DualProblem} we derive dual representations of~$\Pbt$
and their optimality conditions. These are used in
section~\ref{VariationalProperties} to characterize the variational
properties of the value function $v$.  The conjugate, horizon, and
perspective functions arise naturally as part of the analysis, and we
present a calculus (section~\ref{sec:preliminaries}) for these functions
that allows explicit computation of the subdifferential of $v$ for
large classes of misfit functions $\rho$ and regularization functions
$\phi$ (see section~\ref{2 families}).

One of the motivating problems for the general analysis and methods we
present is the treatment of a robust misfit function $\rho$ (such as
the popular Huber penalty) in the context of sparsity promotion, which
typically involves a nonsmooth regularizer $\phi$.  In
section~\ref{Numerics} we demonstrate that the sensitivity analysis
can be applied to solve a sparse nonnegative denoising problem with
convex and nonconvex robust misfit measures.

The proofs of all of the results are relegated to the appendix
(section~\ref{sec:proofs-results}).

\subsection{Notation}

For a matrix $A\in\R^{m\times n}$, the image and inverse image of the
sets $E$ and $F$, respectively, are given by the sets
\[
  AE = \set{ y }{ y = Ax,\, x\in E }
  \quad\mbox{and}\quad
  A^{-1}F = \set{ x }{Ax \in F }.
\]
For a function $\map{p}{\Rn}{\eR}$, its epigraph is denoted $\epi{p}=\set{(x,\mu)}{p(x)\le\mu}$, and its level
set is denoted $\lev{p}{\tau}=\set{x}{p(x)\le\tau}$. 
The function $p$ is said to be proper if $\dom{p}\ne\emptyset$ and closed if $\epi{p}$ is a closed set.
The function
$\indicator{x}{X}$ is the indicator to a set $X$, i.e., $\indicator{x}{X}=0$ if $x\in X$ and
$\indicator{x}{X}=+\infty$ if $x\notin X$.

\section{An inverse function theorem for optimal value functions} \label{IFT}

Let $\map{\psi_i}{X\subseteq\Rn}{\eR}$, $i \in \{1,2\}$, be
arbitrary scalar-valued functions, and
consider the following pair of related problems,
and their associated value functions:
\begin{align}
  \label{P12}\tag*{$\cP_{1,2}(\tau)$}
  v_{1}(\tau) &:=
  \displaystyle{\inf_{x\in X}}\ \psi_1(x) + \indicator{(x, \tau)}{\epi{\psi_2}},
\\\label{P21}\tag*{$\cP_{2,1}(\sigma)$}
    v_{2}(\sigma) &:=
    \displaystyle{\inf_{x\in X}}\ \psi_2(x) + \indicator{(x, \sigma)}{\epi{\psi_1}}.
\end{align}
This pair corresponds to the problems $\Pbt$ and $\Pbs$, defined in
section~\ref{sec:introduction}, with the identifications
\[
 \psi_1(x) = \rho(b-Ax) \text{and} \psi_2(x) = \phi(x).
\]
Our goal in this section is to establish general conditions under
which the value functions $v_1$ and $v_2$ satisfy the inverse-function
relationship
\[
 v_1 \circ v_2 = \mathrm{id},
\]
and for which the %
the pair of
problems~\ref{P12} and~\ref{P21} have the same solution sets. 
The pair of problems~\Pbt
and~\Pbs always satisfy the conditions of the next
theorem, which applies to functions that are not necessarily convex.

\begin{theorem}\label{IVT}
  Let $\map{\psi_i}{X\subseteq\Rn}{\eR}$, $i \in \{1,2\}$, be
  as defined in $\cP_{1,2}(\tau)$, and define
  \[
  \cS_{1,2}:=\set{\tau\in\eR}{\emptyset\ne \argmin\cP_{1,2}(\tau)
  \subset\set{x\in X}{\psi_2(x)=\tau}}.
  \] 
  Let $\cS_{2,1}$ be defined symmetrically to $\cS_{1,2}$ by
  interchanging the roles of the indices.  Then, for every $\tau \in
  \cS_{1,2}$,
  \smallskip
  \begin{enumerate}
  \item[\rm(a)] $v_2(v_1(\tau))=\tau$, and
  \item[\rm(b)] $\argmin \cP_{1,2}(\tau)= \argmin\cP_{2,1}(v_1(\tau))
    \subset\Set{x\in X | \psi_1(x)=v_1(\tau)}$.
  \end{enumerate}
  \smallskip
  Moreover, $\cS_{2,1}=\set{v_1(\tau)}{\tau\in\cS_{1,2}}$, and so
  \[
  \set{\big(\tau,v_1(\tau)\big)}{\tau\in \cS_{1,2}}
  \ =\
  \set{\big(v_2(\sigma),\sigma\big)}{\sigma\in\cS_{2,1}}.
  \]
\end{theorem}

\section{Convex analysis}\label{sec:preliminaries}

In order to present the duality results of section~\ref{DualProblem},
we require a few basic tools from convex analysis. There are many
excellent references for the necessary background material, 
with several appearing within the past 10 years. In this
study we make use of Rockafellar~\cite{RTR} and Rockafellar and
Wets~\cite{RTRW}, although similar results can be found elsewhere
\cite{BauCom:2011,Borwein00, Boyd04,EkeTem:1976,Hiriart-Urruty01,Zal:2002}. 
We review the necessary
results here.

\subsection{Functional operations}\label{sec:generated-functions}

The proper convex function $\map{h}{\Rn}{\eR}$ generates the following convex
functions:
\begin{enumerate}
\item \textit{Legendre-Fenchel conjugate} of $h$:
\[
  h^*(y):=\sup_{x}\left[\ip{y}{x}-h(x)\right].
\]
\item \textit{Horizon function} of $h$:
\[
  h^\infty(z):=\sup_{x\in\dom{h}}\left[h(x+z)-h(x)\right].
\]
\item \textit{Perspective function} of $h$:
\[
\tilde h(z,\lam):=
\begin{cases}
   \lam h(\lam^{-1}z) & \mbox{if}\quad \lam>0,
 \\
 \indicator{z}{0} & \mbox{if} \quad \lam = 0,
 \\
 +\infty        & \mbox{if}\quad \lam < 0.
\end{cases}
\]
\item \textit{Closure of the perspective function} of $h$:
\[
h^{\pi}(z,\lam):=
\begin{cases}
   \lam h(\lam^{-1}z) & \mbox{if}\quad \lam>0,
 \\h^\infty(z)        & \mbox{if}\quad \lam=0,
 \\+\infty            & \mbox{if}\quad \lam<0.
\end{cases}
\]
\end{enumerate}

Each of these functions can also be defined
by considering the epigraphical perspective and properties of convex
sets. Indeed, the horizon function $h^\infty$ is usually defined to be the
function whose epigraph is the horizon cone of the epigraph of $h$
(see section \ref{cones} below). 
The definition given above is a consequence of \cite[Theorem 8.5]{RTR}.

The perspective function of $h$ is the
positively homogeneous function generated by the convex function $\hat
h(x,\lam):=h(x)+\indicator{\lam}{\{1\}}$~\cite[pp.~35 and~67]{RTR}.
If $h$ is additionally closed and proper, then so are $h^*$ (Theorem
12.2), $h^\infty$ (Theorem 8.5), and $h^{\pi}$ (Corollary
8.5.2), where these results are from
Rockafellar~\cite{RTR}. 

Note that for every closed, proper and convex function $h$, the associated
horizon and perspective function, $h^\infty$ and $h^\pi$, are
positively homogeneous and so can be represented as the support
functional for some convex set \cite[Theorem 13.2]{RTR}.  Moreover, if
$h$ is a support function, then $h^\infty=h^\pi=h$.

\subsection{Cones}\label{cones}

We associate the following cones with a convex set $C$ and a convex function $h$.
\begin{enumerate}
\item{\it Polar cone:} The polar cone of $C$ is denoted by 
\[
C^\circ := \set{x^* }{\ip{x^*}{x} \leq 0\;\ \forall x \in C }.
\]
\item{\it Recession cone:} The recession cone of $C$ is denoted by
\[
C^\infty := \set{x}{ C + x \subset C} = \set{x}{y + \lambda x \in C\;\ \forall \lambda \geq 0, \forall y \in C}.
\]
\item{\it Barrier cone:} The barrier cone of $C$ is denoted by
\[
\barrier{C} := \set{x^*}{\mbox{for some $\beta\in \bR, \ip{x}{x^*} \leq \beta\;\ \forall x \in C$}}.
\]
\item{\it Horizon cone of $h$:} The horizon cone~\cite[Theorem 8.7]{RTR} of $h$ is denoted by
\[
\hzn{h}:=\set{y}{h^\infty(y)\le 0}=[\lev{h}{\tau}]^\infty\quad\forall\ \tau>\inf h.
\]
\end{enumerate}

\noindent
A further excellent reference for horizon cones and functions is
\cite{AusTeb:2003} where they are referred to as {\it asymptotic} cones and functions.

\subsection{Calculus rules} \label{sec:calculus}

The conjugate, horizon, and perspective transformations defined in
section~\ref{sec:generated-functions} posses a rich calculus. We use
this calculus to obtain explicit expressions for the functions
$\rho^*,\ \phi^*, (\phi^*)^\infty$ and $(\phi^*)^\pi$, which play a
crucial role in the applications of section~\ref{2 families}.  The
calculus for conjugates and horizons is developed in many references
(e.g., \cite{BauCom:2011,Borwein00,
  Boyd04,EkeTem:1976,Hiriart-Urruty01,Zal:2002}); specific citations
from~\cite{RTR} are provided.  In order to establish the perspective
calculus rules for affine composition and the inverse linear image, we
note that addition is a special case of affine composition, and that
infimal convolution is a special case of inverse linear image.  Hence,
we need only establish the perspective calculus formulas for affine
composition and the inverse linear image:  the formula for affine-composition
follows from \cite[Theorem 9.5]{RTR} and the definition of the
perspective transformation; the formula for inverse linear image is
established in section~\ref{app:inverse linear image}.

\paragraph{\bf Affine composition} Let $\map{p}{\Rm}{\eR}$ be a closed,
proper and convex function, $A\in\Rmn$, and $b\in\Rm$, such that
$(\Ran{A}-b)\cap\ri{\dom{p}}\ne\emptyset$. Let
\begin{flalign*}
   &&h(x)&:=p(Ax-b). &
   \intertext{Then}
   && h^*(y) &= \inf_{A\T u=y}[\ip{b}{u}+p^*(u)], &\llap{\cite[Theorem 16.3]{RTR}}
\\ && h^\infty(z) &= p^\infty(Az), &\llap{\cite[Theorem 9.5]{RTR}}
\\ && h^{\pi}(x,\lam) &= p^{\pi}(Ax-\lam b,\lam),
\end{flalign*}
where, for $\lam=0$,
\[
  h^{\pi}(x,0) =p^{\pi}(Ax,0)=p^\infty(Ax).
\]
All three functions are closed, proper and convex. The derivation of $h^*$ also makes
use of the observation that
\begin{equation}\label{conjugate shift}
\mbox{if} \quad g(x):=h(x-b), \text{then} g^*(v)=h^*(v)+\ip{v}{b}.
\end{equation}

\paragraph{\bf Inverse linear image} 
Let $\map{p}{\Rn}{\eR}$ be
closed, proper and convex, and let $A\in\Rmn$. Let
\[
h(w):=\inf_{Ax=w}p(x)
\]
be the inverse linear image of $p$ under $A$. Then
\begin{flalign*}
   &&  h^*(y)&=p^*(A^Ty). &\llap{\cite[Theorem 16.3]{RTR}}
   \intertext{If $(A^T)^{-1}\ri{\dom{p^*}}\ne\emptyset$, then}
   && h^\infty(z)&=\inf_{Ax=z}p^\infty(x), &\llap{\cite[Theorem 9.2]{RTR}}
\\ && h^{\pi}(w,\lam)&=\inf_{Ax=w}p^{\pi}(x,\lam),&\llap{(Proof in section~\ref{app:inverse linear image})}
\end{flalign*}
where all of the functions $h,\ h^*,\ h^\infty,$ and $h^\pi$ are closed, proper and convex.
\smallskip

\paragraph{\bf Addition} Let $\map{h_i}{\Rn}{\eR}$, for $i=1,\dots,m$, be
closed proper convex functions. If $h:=h_1+\dots+h_m$ is not
identically $+\infty$, then
\begin{flalign*}
   && h^\infty&=h_1^\infty+\dots+h_m^\infty, & \llap{\cite[Theorem 9.2]{RTR}}
\\ && h^{\pi}&=h_1^{\pi}+\dots+h_m^{\pi}, &
\end{flalign*}
where both are closed, proper and convex.
Moreover, if $\bigcap_{i=1}^m\ri{\dom{h_i}}\ne\emptyset$, then
\begin{flalign*}
&& h^*&=h^*_1\infc \cdots\infc h^*_m & \llap{\cite[Theorem 16.4]{RTR}}
\end{flalign*}
is closed, proper and convex, where $\infc$ denotes
infimal convolution.

\paragraph{\bf Infimal convolution}
Let $\map{h_i}{\Rn}{\eR}$, for $i=1,\dots,m$, be closed, proper and convex
functions. Let $h:=h_1\infc \cdots \infc h_m$.  Then $h^*=h_1^*+ \dots +
h^*_m$, and
\begin{flalign*}
 && \mbox{if}\ \ \bigcap_{i=1}^m\ri{\dom{h_i^*}}\ne\emptyset, \text{then}
  h^\infty&=h_1^\infty\infc \cdots \infc h^\infty_m, \qquad\qquad\qquad
  &\llap{\cite[Corollary 9.2.1]{RTR}}
\end{flalign*}
and
\[
h^{\pi}(x,\lam)=\inf_{\sum_{i=1}^mx_i=x}
\left[h_1^{\pi}(x_1,\lam)+\dots+h_m^{\pi}(x_m,\lam)\right].
\]
All three functions are closed, proper and convex.

\section{The dual problem}
\label{DualProblem}

For our analysis, it is convenient to consider the (equivalent) epigraphical formulation
\begin{equation}
  \tag{$\cP$}\label{ConstrainedClass}
  v(b,\tau) = \minim_x\ f(x,b,\tau)
\end{equation}
of \Pbt, where
\begin{equation*}\label{perturbation family}
  f(x,b,\tau):=\rho(b - Ax) + \indicator{(x,\tau)}{\epi{\phi}}.
\end{equation*}
Because the functions $\rho$ and $\phi$ are convex, it immediately
follows that $f$ is also convex. This fact gives the convexity of the
value function $v$, since it is the inf-projection of the objective
function in $x$~\cite[Theorem~5.3]{RTR}.

We use a duality framework derived from the one described in Rockafellar and
Wets~\cite[Chapter 11, Section H]{RTRW}, and
associate with~\ref{ConstrainedClass} %
its dual problem and corresponding dual value function:
\begin{equation}\label{dual1}
  \tag{$\cD$}
  \hv(b,\tau) := \maxim_{u,\mu}\ \ip{b}{u}+\tau\mu-f^*(0,u,\mu).
\end{equation}
To derive this dual from~\cite[Theorem 11.39]{RTRW},
define 
\[
f_{(b,\tau)}(x,\Delta b,\Delta\tau):=f(x,b+\Delta b,\tau+\Delta \tau)\ .
\]
Then, by \eqref{conjugate shift}, 
$f^*_{(b,\tau)}(v,u,\mu)=f^*(v,u,\mu)-\ip{b}{u}-\tau\mu$. Substituting this expression 
into~\cite[Theorem 11.39]{RTRW} gives $\cD$.

The dual $\cD$ is the key to understanding the variational behavior of
the value function.  To access these results we must compute the
conjugate of $f$. For this it is useful to have an alternative
representation for the support function of the epigraph, which is the
conjugate of the indicator function appearing in
$f$. %

\subsection{Reduced dual problem}
In Theorem~\ref{reduced dual theorem}, we derive an
equivalent representation of the dual problem~\ref{dual1} in terms of
$u$ alone. This is the \emph{reduced} dual problem
for~\ref{ConstrainedClass}. We first present a result about conjugates
for epigraphs and lower level sets.

\begin{lemma}[Conjugates for epigraphs and lower level sets]
\label{lem:conjugate of epi indicator}
Let $\map{h}{\Rn}{\eR}$ be closed, proper and convex. Then
\begin{subequations}
\begin{align}
  \support{(y,\mu)}{\epi{h}}&=(h^*)^{\pi}(y,-\mu) \label{epi support and perspective},
\\\support{y}{\lev{h}{\tau}}&=\cl{\inf_{\mu\ge 0}\left[\tau\mu+(h^*)^\pi(y,\mu)\right]}.
  \label{level support}
\end{align}
\end{subequations}
\end{lemma}

Expressions~\eqref{level support} and~\eqref{epi support and
  perspective} are easily derived from the case where $\tau=0$ which
is established in~\cite[Theorem 13.5]{RTR} and~\cite[Corollary
13.5.1]{RTR}, respectively. In \cite{RTR}, it is shown that~\eqref{epi
  support and perspective} is a consequence of~\eqref{level
  support}. In section~\ref{sec:proofs-results} we provide a different
proof of Lemma~\ref{lem:conjugate of epi indicator} where it is shown
that~\eqref{level support} follows from~\eqref{epi support and
  perspective}. The arguments provided in the proof are instructive
for later computations.

The conjugate $f^*(y,u,\mu)$ of the perturbation function
$f(x,b,\tau)$ defined in~\ref{ConstrainedClass} is now easily computed:
\begin{align}
f^*(y,u,\mu)
   &=\sup_{x,b,\tau}\left[\ip{y}{x}+\ip{u}{b}+\mu\tau -\rho(b-Ax)-\indicator{(x,\tau)}{\epi{\phi}}\right]
\nonumber\\ 
&=\sup_{x,w,\tau}\left[\ip{y}{x}+\ip{u}{w+Ax}+\mu\tau -\rho(w)-\indicator{(x,\tau)}{\epi{\phi}}\right]
\nonumber\\ 
&=\sup_{x,\tau}\left[\ip{y+A^Tu}{x}+\mu\tau-\indicator{(x,\tau)}{\epi{\phi}}\right]
+\sup_{w}\left[\ip{u}{w}-\rho(w)\right]
\nonumber\\ 
&=(\phi^*)^{\pi}(y+A^Tu,-\mu)+\rho^*(u),\label{f conj}
\end{align}
where the final equality follows from~\eqref{epi support and
  perspective}.  With this representation of the conjugate of $f$, we
obtain the following equivalent representations for the dual
problem~\ref{dual1}. The representation labeled $\cD_r$ is of particular importance to our
discussion.  We refer to $\cD_r$ as the {\it reduced dual}.

\begin{theorem}[Dual representations]\label{reduced dual theorem}
For problem~\ref{ConstrainedClass} define the functions
\begin{align*}
  g_\tau(u)&:=\rho^*(u)+\support{A^Tu}{\lev{\phi}{\tau}}
\\p_\tau(s,\mu)&:=\tau\mu + (\phi^*)^\pi(s,\mu).
\end{align*}
Then the value function for \ref{dual1} has the following
equivalent characterizations:
\begin{subequations}
\begin{align}
\hv(b,\tau)&=
\sup_{u}\left[\ip{b}{u}-\rho^*(u)-\inf_{\mu\ge 0}\ p_\tau(A^T u,\mu)\right]\nonumber
\\ &= 
\sup_{u}\left[\ip{b}{u}-\rho^*(u)-\support{A^Tu}{\lev{\phi}{\tau}}\right]\label{GeneralValueDualr}\tag{$\cD_r$}
\\ &= g_\tau^*(b)\label{g conj}
\\ &=  \cl{v(\cdot,\tau)}(b),\label{g conj2}
\end{align}
\end{subequations}
where the closure operation in the~\eqref{g conj2} refers to the
lower semi-continuous hull of the convex function $b\mapsto
v(b,\tau)$.  In particular, this implies the weak duality inequality
$\hv(b,\tau)\le v(b,\tau)$. Moreover, if the function $\rho$ is
differentiable, the solution $u$ to~\ref{GeneralValueDualr} is
unique.
\end{theorem}

In the large-scale setting, the primal problem \Pbt is usually solved
using a primal method
that does not give direct access to the multiplier $\bmu$
for the inequality constraint $\phi(x)\le\tau$.  For example, \Pbt may
be solved using a variant of the gradient projection algorithm.
However, one can still obtain an approximation to the optimal dual
variable $\bu$ in~\ref{GeneralValueDualr}, typically through the
residual corresponding to the current iterate. For this reason, one needs a way
to obtain an approximation to $\bmu$ from an approximation to
$\bu$ (i.e., given $\bu$, compute $\bmu$).  Lemma \ref{lem:conjugate
  of epi indicator} and Theorem~\ref{reduced dual theorem} show that
this can be done by solving the problem $\inf_{\mu\ge 0}\ p_\tau(A^T
\bu,\mu)$ for $\bmu$.  Indeed, in the sequel we show that in many
important cases there is a closed form expression for the solution
$\bmu$.  The following lemma serves to establish a precise
relationship between the solution $\bu$ to the reduced
dual~\ref{GeneralValueDualr} and the solution pair $(\bu,\bmu)$ to the
dual \ref{dual1}.

\begin{lemma}\label{mu existence}
  Let $\phi$ be as in~\ref{ConstrainedClass} with $\tau > \inf \phi$ and $\bx\in
  \lev{\phi}{\tau}$.

\medskip

\noindent 1. \label{ItemOne} For every $s$, we have
\begin{equation}
\label{LevelSupportInf}
\support{s}{\lev{\phi}{\tau}}\leq\inf_{\mu\ge 0}p_\tau(s,\mu).
\end{equation}

\medskip

\noindent 2.
Let $(\bx, \bs)$ satisfy 
  $\bs\in\ncone{\bx}{\lev{\phi}{\tau}}$
  and define %
\begin{equation*}
  S_1=\argmin_{\mu\ge 0}\, p_\tau(\bs,\mu) \text{and}
  S_2= \set{\bar\mu\ge 0}{\begin{aligned}
\bs&\in\bmu^+\sd\phi(\bx)\\
0 &= \bmu(\phi(\bx) - \tau) \end{aligned}},
\end{equation*}
where, for $x\in \dom{\phi}$, 
\[
\mu^+\partial \phi(x):=\begin{cases}
\set{\mu z}{z\in\partial \phi(x)} ,&\mbox{if $\mu >0$ and $x \in \dom{\sd{\phi}}$},\\
\ncone{x}{\dom{\phi}} ,&\mbox{if $\mu=0$ or $\sd{\phi}(x)=\emptyset$}.
\end{cases}
\]
If either %
$S_1$ or $S_2$ is non-empty, then $S_1=S_2$
and equality holds in \eqref{LevelSupportInf}.
\end{lemma}

In \Zal the object $\mu^+\sd \phi(x)$ is denoted as $\sd(\mu \phi)(x)$ \cite[page 141]{Zal:2002}, where 
\[
(\mu \phi)(x):=\begin{cases}
\mu \phi(x),&\mbox{ if $\lam >0$, and}\\
\indicator{x}{\dom{\phi}},&\mbox{ if $\lam=0$.}
\end{cases}
\]
We choose the notation $\mu^+\sd\phi(x)$ to emphasize that there is an underlying limiting
operation at play, e.g. see \cite[Definition 8.3 and Proposition 8.12]{RTRW}.

The final lemma of this section concerns conditions under which solutions
to~\ref{ConstrainedClass} and~\ref{GeneralValueDualr} exist. This is closely tied
to the horizon behavior of these problems and the notion of coercivity.

\begin{defn}
A function $\map{h}{\Rn}{\eR}$ is said to be $\alf$-coercive if
\[
\lim_{\norm{x}\rightarrow\infty}\frac{f(x)}{\norm{x}^\alf}=+\infty.
\]
In particular, $h$ is said to be $0$-coercive, or simply coercive, if
$\lim_{\norm{x}\rightarrow\infty}f(x)=\infty$.
\end{defn}

\begin{lemma}[Coercivity of primal and dual objectives]\label{coercivity 1}
\begin{subequations}
\smallskip

\noindent 1.
The objective function $f(\cdot,b,\tau)$
  of~\ref{ConstrainedClass} is coercive if and only if
  \begin{equation}
    \label{eq:primal coercivity}
    \hzn{\phi}\cap[-A^{-1}\hzn{\rho}]=\{0\}.
  \end{equation}

\smallskip

\noindent 2.
The objective function of the reduced
  dual~\ref{GeneralValueDualr} is coercive if and only if
\begin{equation}
  \label{eq:dual coercivity}
  b\in\intr{\dom{\rho}+A\lev{\phi}{\tau}}.
\end{equation}
\end{subequations}
\end{lemma}

\section{Variational properties of the value function}\label{VariationalProperties}

Using \ref{dual1} and representation of the conjugate of the
objective of~\ref{ConstrainedClass} (cf.~\eqref{f conj}), we can
specialize \cite[Theorem 11.39]{RTRW} to obtain a characterization of
the subdifferential of the value function, as well as sufficient
conditions for strong duality.

\begin{theorem}[Strong duality and subgradient of the value
  function] \label{value function sd} Let $v$ and $\hv$ be as
  in~\ref{ConstrainedClass} and \ref{dual1}, respectively.  It is
  always the case that
  \begin{align*}
  v(b,\tau)&\ge\hv(b,\tau) \qquad\mbox{(weak duality).}
  \intertext{If $(b,\tau)\in\intr{\dom{v}}$, then}
  v(b,\tau)&=\hv(b,\tau) \qquad\mbox{(strong duality)}
  \end{align*}
and
  $\partial v(b,\tau)\ne\emptyset$ with
  \begin{equation*}
  \partial v(b,\tau)
  :=\argmax_{u,\ \mu\ge0}\ [\ip{b}{u}-\rho^*(u) - p_\tau(A^T u,-\mu)]\,.
  \end{equation*}
  Furthermore, for fixed $(b,\tau)\in\Rm\times\R$,
  \[ \dom{f(\cdot,b,\tau)}\neq \emptyset \iff b \in
  \dom{\rho}+A(\lev{\phi}{\tau}).  \]
  In particular, this implies that
  \[
  (b,\tau)\in\intr{\dom{v}} \iff b\in\intr{\dom{\rho}+A(\lev{\phi}{\tau})}.
  \]
\end{theorem}

We now derive a characterization of the subdifferential $\partial
v(b,\tau)$ based on the solutions of the reduced dual
\ref{GeneralValueDualr}.

\begin{theorem}[Value function
  subdifferential]\label{reduced primal-dual}
\begin{subequations} \label{eq:reduced primal-dual}
  Suppose that
\begin{gather}
  b\in\ri{\dom{\rho}}+A\,\ri{\lev{\phi}{\tau}} \quad\mbox{and}\label{primal regularity}
\\\ri{\dom{\rho^*}}\cap[A^{-T}\ri{\barr{\lev{\phi}{\tau}}}]\ne\emptyset.
\label{dual regularity}
\end{gather}

\smallskip

\noindent 1. If the pair $(\bx,\bu)$ satisfies
\begin{equation}\label{ReducedOptCond}
\bx\in\lev{\phi}{\tau},\ \bu\in\sd \rho(b-A\bx)\mbox{ and } A^T\bu\in\ncone{\bx}{\lev{\phi}{\tau}},
\end{equation}
then $\bx$ solves~\ref{ConstrainedClass}
and $\bu$ solves \ref{GeneralValueDualr}.

\smallskip

\noindent 2. If $\bx$ solves~\ref{ConstrainedClass} and \eqref{primal regularity} holds,
there exists
$\bu$ such that $(\bx,\bu)$ satisfies~\eqref{ReducedOptCond}.

\smallskip

\noindent 3.
If $\bu$ solves \ref{GeneralValueDualr} and \eqref{dual regularity} holds,
there exists $\bx$ such that  $(\bx,\bu)$ satisfies~\eqref{ReducedOptCond}.

\smallskip

\noindent 4. If either \eqref{eq:primal coercivity} and \eqref{primal
  regularity} holds, or \eqref{eq:dual coercivity} and \eqref{dual
  regularity} holds, then $\sd v(b,\tau)\ne\emptyset$ and
$\argmin_{\mu\ge 0}\ p_\tau(A^T \bu,\mu)\,\ne\,\emptyset$ for all
$(\bx,\bu)\in\Rn\times\Rm$ satisfying \eqref{ReducedOptCond} with
\begin{align}
\label{subdiff explicit}
\sd v(b,\tau)&=
  \set{ \begin{pmatrix}\hfill\bu\\-\bmu\end{pmatrix}}{
  \begin{array}{c}(\bx,\bu)\!\in\!\Rn\!\times\!\Rm \,\mbox{satisfy \eqref{ReducedOptCond} and}\\
    \bmu \in \argmin_{\mu\ge 0}\ p_\tau(A^T \bu,\mu)
\end{array}
}\quad
\\
\label{subdiff explicit 2}
&=
  \set{
    \begin{pmatrix}\hfill\bu\\-\bmu\end{pmatrix}}{
    \begin{array}{c}
    \exists\, \bx\in\lev{\phi}{\tau}\, \mbox{s.t.}\
    0\!\in\! -A^T\bu\!+\!\bmu^+\partial\phi(\bx)\\ \mbox{ where }
    \bu\in\sd\rho(b-A\bx),\\ \bmu\ge 0,\mbox{ and }\bmu(\phi(\bx)-\tau)=0
\end{array}
}.
\end{align}
\end{subequations}
\end{theorem}

The representation \eqref{subdiff explicit 2} expresses the elements
of $\partial v(b,\tau)$ in terms of classical Lagrange multipliers
when $\bmu>0$, and extends the classical theory when $\bmu=0$. (See
Lemma~\ref{mu existence} for the definition of
$\mu^+\partial\phi(x)$.)
Because $v$ is convex, it is subdifferentially regular, and so for
fixed $b$, we can obtain the subdifferential of $v$ with respect to
$\tau$ alone \cite[Corollary 10.11]{RTRW}, i.e.,
\[
\partial_\tau v(b,\tau) =
\set{\omega}{\begin{pmatrix}u\\\omega\end{pmatrix}\in\partial v(b,\tau)}.
\]

\section{Applications}\label{2 families}

In this section we apply the calculus rules of
section~\ref{sec:calculus} in conjunction with Theorem~\ref{reduced
  primal-dual} to evaluate the subdifferential of the value function
in three important special cases: where $\phi$ is a
gauge-plus-indicator (section~\ref{GaugeIndicator}), a quadratic
support function (section \ref{LinearQuadratic}), and an affine
composition with a quadratic support function (section~\ref{EPLQ}). In
all cases we allow $\rho$ to be an arbitrary convex function.

\subsection{Gauge-plus-indicator}
\label{GaugeIndicator}
The case where $\rho$ is a linear least-squares objective and $\phi$
is a gauge function is studied in \cite{BergFriedlander:2011}.  We
generalize this case by allowing the convex function $\rho$ to be a
possibly non-smooth and non-finite-valued, and take
\begin{equation}\label{case 1}
\phi(x):=\gauge{x}{U}+\indicator{x}{X},
\end{equation}
where $U$ is a nonempty closed convex set containing the origin.
Here, $\gauge{x}{U}$ is the gauge function defined in
\eqref{LittleLS}. It is evident from the definition of a gauge that
$\phi$ is also a gauge if and only if $X$ is a convex cone. Since
$0\in U$, it follows from \cite[Theorem 14.5]{RTR} that
$\gauge{\cdot}{U}=\support{\cdot}{U^\circ}$, where
\[
U^\circ=\set{v}{\ip{v}{u}\le 1\ \forall\, u\in U}
\] 
is the polar of the set $U$.

Observe that the requirement $x\in X$ is unaffected by varying
$\tau$ in the constraint $\phi(x)\le\tau$. Indeed, the
problem~\ref{ConstrainedClass} is unchanged if we replace $\rho$ and
$\phi$ by
\begin{equation}\label{case 1A}
\hrho(y,x):=\rho(y)+\indicator{x}{X}\quad\mbox{ and }\quad\hphi(x):=\gauge{x}{U},
\end{equation}
with $A$ and $b$ replaced by
\[
\hb:=\begin{pmatrix} b\\ 0\end{pmatrix}\quad\mbox{ and }\quad
\hA:=\begin{bmatrix}\hfill A\\ -I\end{bmatrix}.
\]
Hence, the generalization of \cite{BergFriedlander:2011} discussed here
only concerns the application to more general convex functions $\rho$.

There are two ways one can proceed with this application. One can use
$\phi$ as given in \eqref{case 1} or use $\hrho$ and $\hphi$ as
defined in \eqref{case 1A}.  We choose the former in order to
highlight the presence of the abstract constraint $x\in X$.  But we
emphasize---regardless of the formulation chosen---the end result
is the same.

\begin{lemma}\label{case 1 details}
Let $\phi$ be as given in \eqref{case 1}. 
The following formulas hold:
\begin{subequations}
\begin{eqnarray}
\gauge{\cdot}{U}&=&\support{\cdot}{U^\circ},\label{gauge=support}\\
\dom{\gauge{\cdot}{U}}&=&\cone{U}=\barrier{U^\circ},\label{gauge dom}\\
\dom{\phi}&=&\cone{U}\cap X,\label{c1 dom phi}\\
\lev{\phi}{\tau}&=&(\tau U)\cap X,\label{c1 lev phi}\\
\hzn{\phi}&=&U^\infty\cap X^\infty,\mbox{ and }\label{c1 hzn phi}\\
\cl{\barrier{\lev{\phi}{\tau}}}&=&\cl{\barrier{U}+\barrier{X}}.\label{c1 bar lev phi}
\end{eqnarray}
\end{subequations}
If it is further assumed that %
\begin{equation}\label{c1 cq}
\ri{\tau U}\cap\ri{X}\ne\emptyset,
\end{equation}
then we also have
\begin{subequations}
\begin{eqnarray}
 \phi^*(z)&=&{\min_s[\support{z-s}{X}+\indicator{s}{U^\circ}]},\label{c1 conj}\\
(\phi^*)^\pi(z,\mu)&=&{\min_s[\support{z-s}{X}+\indicator{s}{\mu U^\circ}]},
\label{c1 conj pi}\\
\support{z}{\lev{\phi}{\tau}}&=&{\min_{\mu\ge 0}[\tau\mu +(\phi^*)^\pi(z,\mu)]}
\label{c1 support lev phi a}\\
&=&{\min_s[\support{z-s}{X}+\tau\gauge{s}{U^\circ}]},
\mbox{ and}\qquad\quad \label{c1 support lev phi}\\
\ncone{x}{\lev{\phi}{\tau}}&=&\ncone{x}{\tau U}+\ncone{x}{X}.\label{c1 ncone lev phi}
\end{eqnarray}
\end{subequations}
If $\bs$ minimizes \eqref{c1 support lev phi}, then
$\bmu:=\gauge{\bs}{U^\circ}$ minimizes \eqref{c1 support lev phi a}.
\end{lemma}

By Theorem \ref{value function sd}, the subdifferential of $v(b,\tau)$
is obtained by solving the dual problem~\eqref{GeneralValueDual} or
the reduced dual \ref{GeneralValueDualr}.  When $\phi$ is given by
\eqref{case 1}, the results of Lemma~\ref{case 1 details} show that
the dual and the reduced dual take the form
\begin{align}
\label{case 1 dual}
\sup_{u,\mu}&\left[\ip{b}{u}+\tau\mu- (\phi^*)^{\pi}(A^Tu,-\mu)-\rho^*(u)\right]
\\ &= \nonumber
\sup_{u}\left[\ip{b}{u}-\rho^*(u)-\support{A^Tu}{\lev{\phi}{\tau}}\right]
\\ &= \nonumber
\sup_{u}\left[\ip{b}{u}-\rho^*(u)-{\min_s[\support{A^Tu-s}{X}+\tau\gauge{s}{U^\circ}]}\right]
\\ &=
\sup_{u,s}\left[\ip{b}{u}-\rho^*(u)-\support{A^Tu-s}{X}-\support{s}{\tau U}\right].\label{case 1 reduced}
\end{align}
Moreover, if $(\bu,\bs)$ solves \eqref{case 1 reduced}, then $(\bu,\bmu)$ solves \eqref{case 1 dual}
with $\bmu=-\gauge{\bs}{U^\circ}$, and
\[
 (\bu,-\gauge{\bs}{U^\circ})\in\sd v(b,\tau).
\]
We have the following version of Theorem \ref{reduced primal-dual}
when $\phi$ is given by \eqref{case 1}.

\begin{theorem}\label{case 1 primal-dual}
Let $\phi$ be given by \eqref{case 1} under the assumption that
\eqref{c1 cq} holds, and consider the following two conditions:
\begin{equation}\label{c1 primal regularity}
b\in\ri{\dom{\rho}+A[{\tau U}\cap{X}]}=\ri{\dom{\rho}}+A[\ri{\tau U}\cap\ri{X}]
\end{equation}
and
\begin{equation}\label{c1 dual regularity}
\exists\ \hu\in\ri{\dom{\rho^*}}
\quad\mbox{such that}\quad
A^T \hu\in\ri{\barrier{U}}+\ri{\barrier{X}}.
\end{equation}

\smallskip

\noindent 1. \label{case 1.1}
If the triple $(\bx,\bu,\bs)$ satisfies
\begin{subequations} \label{cases1ab}
\begin{gather}
\bu\in\sd \rho(b-A\bx),\qquad  \bx\in X\cap(\tau U),
\\
\bs\in\ncone{\bx}{\tau U},\quad\mbox{and}\quad
A^T\bu-\bs\in\ncone{ \bx}{X},
\end{gather}
\end{subequations}
then $\bx$ solves $\Pbt$ %
and $(\bu,\bs)$ solves \eqref{case 1 reduced}.

\smallskip

\noindent 2. \label{case 1.2}
If $\bx$ solves $\Pbt$ %
and \eqref{c1 primal regularity} holds,
then there exists
a pair $(\bu,\bs)$ such that $(\bx,\bu,\bs)$ satisfies~\eqref{cases1ab}.

\smallskip

\noindent 3. \label{case 1.3}
If $(\bu,\bs)$ solves \eqref{case 1 reduced} and \eqref{c1 dual regularity} holds,
then there exists $\bx$ such that  $(\bx,\bu,\bs)$ satisfies~\eqref{cases1ab}.

\smallskip

\noindent 4.
If
either
\begin{equation}\label{c1 coercive primal}
U^\infty\cap X^\infty\cap [-A^{-1}\hzn{\rho}]=\{0\}\quad
\mbox{and \eqref{c1 primal regularity} holds,}
\end{equation}
or
\begin{equation}\label{c1 coercive dual}
b\in\intr{\dom{\rho}+A[{\tau U}\cap{X}]}
\mbox{and \eqref{c1 dual regularity} holds,}
\end{equation}
then $\sd v(b,\tau)\ne\emptyset$ and is given by
\begin{align}\label{c1 subdiff explicit}
\sd v(b,\tau)&=
  \set{
    \begin{pmatrix}\bu\\-\gauge{\bs}{U^\circ}\end{pmatrix}}{
    (\bx,\bu,\bs)\in\Rn\times\Rm\times\Rn \mbox{
    satisfy~\eqref{cases1ab}}}\!\!\quad
    \\ \ \nonumber\\ \label{c1 subdiff explicit 2}
  &=
  \set{
    \begin{pmatrix}\hfill\bu\\-\bmu\end{pmatrix}}{
       \begin{array}{c}
       \exists\ \bx\in X\!\cap\!(\tau U)\mbox{ s.t. }\\
       0\in -A^T\bu+\ncone{\bx}{X}+\bmu^+\sd \gauge{\bx}{U}\mbox{ where}\!\!\\
       \bu\in\sd \rho(b-A\bx),\ 0\le\bmu\mbox{ and }\bmu(\gauge{\bx}{U}-\tau)=0\
        \end{array}
      }. \nonumber
 \end{align}

\end{theorem}

\subsubsection{Gauge penalties}
\label{NormIndicatorSection}
In \cite{BergFriedlander:2011}, the authors study the case where
$\rho$ is a linear least-squares objective, $\phi$ is a gauge
functional, and $X=\R^n$. In this case, \cite[Lemma
2.1]{BergFriedlander:2011} and \cite[Theorem
2.2(b)]{BergFriedlander:2011} can be deduced from \eqref{case 1
  reduced} and \eqref{c1 subdiff explicit}, respectively.  Another
application is to the case where $\rho$ is finite-valued and smooth,
$\phi$ is a norm, and $X$ is a generalized box.  In this case, all of
the conditions of Theorems \ref{value function sd} and \ref{case 1
  primal-dual} are satisfied, solutions to both
$\Pbt$ %
and \eqref{case 1 reduced} exist, and $v$ is differentiable.
In particular, consider the non-negative 1-norm-constrained inversion, where
\[
\label{NNL1}
\phi(x) = \|x\|_1 + \indicator{x}{\bR^n_+},
\]
and $\rho$ is any differentiable convex function.
The subdifferential characterization given in Theorem~\ref{value
 function sd} can be explicitly computed via Theorem~\ref{case 1
 primal-dual}.  In the notation of~\eqref{case 1},
\[
U = \Set{x | \norm{x}_1 \le 1}=:\mathbb{B}_1,
\]
and $X$ in~\eqref{case 1} is
$\R_+^n$.  Since the function $\rho$ is differentiable, the solution
$\bu$ to the dual~\eqref{case 1 reduced} is
unique~\cite[Theorem~26.3]{RTR}. Therefore, Theorem~\ref{case 1
  primal-dual} gives the existence of a {\it unique} gradient
\[
\nabla_b v(b, \tau) = -A^T\nabla \rho(b - A\bar x),
\]
where $\bar x$ is any solution that achieves the optimal value. 
The derivative with respect to $\tau$ is immediately given by
Theorem~\ref{case 1 primal-dual} as
\begin{equation}
\label{HuberTauDerivative}
\nabla_\tau v(b, \tau) = -\gauge{A^T\nabla \rho(b - A\bar x)}{U^\circ} = 
-\|A^T\nabla \rho(b - A\bar x)\|_\infty.
\end{equation}

Note that~\eqref{HuberTauDerivative} has the same algebraic form when $x$ is unconstrained.  
The non-negativity constraint on $x$ is reflected in the derivative only through 
its effect on the optimal point $\bar x$. 

\subsection{Quadratic support functions} \label{LinearQuadratic}

We now consider the case
\begin{equation}\label{LQ phi}
\phi(x):=\sup_{w\in U}\, [\ip{x}{w}-\half\ip{w}{Bw}]\ ,
\end{equation}
where $U\subset\Rn$ is nonempty, closed, and convex with $0\in U$, and
$B\in\Rnn$ is positive semidefinite. We call this class of functions
quadratic support (\RQP) functions. 
This surprisingly rich and useful class is found in many applications.
A deeper study of its properties and uses 
can be found in \cite{ABP:2013}.
Note that the conjugate of $\phi$ is given by
\begin{equation}\label{case 2 phi conj}
\phi^*(w)=\half\ip{w}{Bw}+\indicator{w}{U}.
\end{equation}
If the set $U$ is polyhedral convex, then the function $\phi$ is called a 
piecewise linear-quadratic (PLQ) penalty function \cite[Example 11.18]{RTRW}.
Since $B$ is positive semidefinite there is a matrix $L\in\R^{n\times
  k}$ such that $B=LL^T$ where $k$ is the rank of $B$. Using $L$, the
calculus rules in section~\ref{sec:calculus} give the following
alternative representation for $\phi$:
\begin{eqnarray}
\phi(x)
   &=&\sup_{w\in U}\, [\ip{w}{x}-\half \|L^Tw\|_2^2-\indicator{w}{U}]\nonumber
\\ &=&\inf_{x_1+x_2=x}\left[\support{x_1}{U}+\inf_{Ls=x_2}\half\tnorm{s}^2\right]\nonumber
\\ &=&\inf_{s}\left[\half\tnorm{s}^2+\support{x-Ls}{U}\right]\nonumber
\\ &=&\inf_{s}\left[\half\tnorm{s}^2+\gauge{x-Ls}{U^\circ}\right],\label{LQ gauge}
\end{eqnarray}
where the final equality follows from~\cite[Theorem 14.5]{RTR} since
$0\in U$.
Note that the function class \eqref{LQ phi} includes all gauge functionals for
sets containing the origin.
By \eqref{case 2 phi conj},  it easily follows that
\[
(\phi^*)^\pi(w,\mu)=
\begin{cases}
  \frac{1}{2\mu}\norm{w}_B^2+\indicator{w}{\mu U} & \mbox{if $\mu>0$,}
\\\indicator{w}{U^\infty\cap\Nul{B}}              & \mbox{if $\mu=0$,}
\\+\infty                                         & \mbox{if $\mu<0$,}
\end{cases}
\]
where $\|\cdot\|_B$
denotes the seminorm induced by $B$, i.e.,
\begin{equation*}
  \label{seminorm}
  \|w\|_B := \sqrt{w^TBw}.
\end{equation*}

The next result catalogues important properties of the function $\phi$
given
in~\eqref{LQ phi}.

\begin{lemma}\label{c2 phi}
  Let $\phi$ be given by \eqref{LQ phi} with $\tau>0$. Then
\begin{eqnarray*}
  \dom{\phi}&=&\cone{U^\circ}+\Ran{B} %
  \quad\mbox{ and }\label{dom LQ phi}\\
  \hzn{\phi}&=&{\cone{U}}^\circ, \label{hzn LQ phi}
\end{eqnarray*}
in particular, $\phi$ is coercive if and only if $0\in\intr{U}$.
Moreover,
\begin{align}
\support{w}{\lev{\phi}{\tau}}
   &=\min_{\lambda\ge 0}\, [ \tau\lambda +(\phi^*)^\pi(w,\lam)]
\label{ind conj inf rep}
\\ &=
\begin{cases}
  \tau\gauge{w}{U} + \frac{\|w\|_B^2}{2\gauge{w}{U}} &
  \mbox{if $\gauge{w}{U}> \|w\|_B/\sqrt{2\tau}$,}
\\ \sqrt{2\tau}\|w\|_B &
   \mbox{if $\gauge{w}{U}\leq \|w\|_B/\sqrt{2\tau}$,}
\end{cases}
\label{indicatorConjugate}
\end{align}
where the minimizing $\lam$ in \eqref{ind conj inf rep} is given by
\begin{equation}
\label{minimizingLambda}
\lambda = \max\left\{\gauge{w}{U},\ \frac{\|w\|_B}{\sqrt{2\tau}}\right\}.
\end{equation}
In particular, the formula \eqref{indicatorConjugate} implies that
\begin{equation*}\label{c2 dom support phi}
\barrier{\lev{\phi}{\tau}}=
\dom{\left(\support{\cdot}{\lev{\phi}{\tau}}\right)}=\dom(\gauge{\cdot}{U})=\cone{U}.
\end{equation*}
\end{lemma}

We now apply Theorem \ref{reduced primal-dual} to the case where
$\phi$ is given by \eqref{LQ phi}.

\begin{theorem}\label{case 2 sd v}
Let
$\phi$ be given by~\eqref{LQ phi}, and consider the following two
conditions:
\begin{equation}\label{c2 primal regularity}
  \exists\ \hx\in\ri{\dom{\phi}} %
  \quad\mbox{such that}\quad \phi(\hx)<\tau\mbox{ and }
  b-A\hx\in\ri{\dom{\rho}}
\end{equation}
and
\begin{equation}\label{c2 dual regularity}
  \exists\ \hu\in\ri{\dom{\rho^*}}
  \quad\mbox{such that}\quad
  A^T\hu\in\ri{\cone{U}}.
\end{equation}

\medskip

\noindent 1. \label{case 2.1} If the pair $(\bx,\bu)$ satisfy
  \begin{equation}\label{c2 opt cond}
 \bx\in\lev{\phi}{\tau},\    \bu\in\sd\rho(b-A\bx)
    \text{ and } 
    A^T\bu\in\ncone{\bx}{\lev{\phi}{\tau}},%
  \end{equation}
then $\bx$ solves~ $\Pbt$ and $\bu$ solves
\ref{GeneralValueDualr}.

\medskip

\noindent 2. \label{case 2.2} If $\bx$ solves~ $\Pbt$ and
  \eqref{c2 primal regularity} holds, then there exists $\bu$ such
  that \eqref{c2 opt cond} holds.

\medskip

\noindent 3. \label{case 2.3} If $\bu$ solves \ref{GeneralValueDualr} and
  \eqref{c2 dual regularity} holds, then there exists $\bx$ such
  that \eqref{c2 opt cond} holds.

\medskip

\noindent 4.
If either
\begin{equation}
\label{c2 primal existence}
\cone{U}^\circ\cap[-A^{-1}\hzn{\rho}]=\{0\}\mbox{ and
\eqref{c2 primal regularity} holds, }
\end{equation}
or
\begin{equation}\label{c2 dual existence}
b\in\intr{\dom{\rho}+A\lev{\phi}{\tau}}\mbox{ and \eqref{c2 dual regularity} holds,}
\end{equation}
then $\sd v(b,\tau)\ne\emptyset$ and is given by
\begin{align}\nonumber%
\sd v(b,\tau)&=\set{\begin{pmatrix}\hfill\bu\\ -\bmu\end{pmatrix} }{
  \begin{array}{c}\exists\ \bx \mbox{ s.t. } (\bx,\bu)\mbox{ satisfy \eqref{c2 opt cond} and  }\\
    \bmu=\max\left\{\gauge{A^T\bu}{U},
      \|A^T\bu\|_B/\sqrt{2\tau}\right\}\end{array}}\quad
\\
&= \set{\begin{pmatrix}\hfill\bu\\ -\bmu\end{pmatrix}}{
\begin{array}{c}
\exists \ \bx\in\lev{\phi}{\tau}\mbox{ s.t. }0\in -A^T\bu+\bmu^+\sd\phi(\bx)
\mbox{ where}\\ \bu\in\sd\rho(b-A\bx),\ 0\le\bmu\mbox{ and }\bmu(\phi(\bx)-\tau)=0
\end{array}}.
\label{c2 sd 2}
\end{align}
\end{theorem}

In the following corollary we exploit the structure of $\phi$ to refine the multiplier description
of the $\sd v(b,\tau)$ given in~\eqref{c2 sd 2}.

\begin{corollary}\label{c2 multiplier cor}
Consider the problem $\Pbt$ with $\phi$ given
by~\eqref{LQ phi}.
A pair $(\bx,\bu)$ satisfies \eqref{c2 opt cond}
if and only if $\bx\in\lev{\phi}{\tau}$, $\bu\in\sd\rho(b-A\bx)$, and either
\begin{subequations} \label{c2 multiplier rule}
\begin{gather} 
  A^T \bu \in \ncone{\bx}{\dom{\phi}},
  \quad\mbox{or}  \label{c2 multiplier rule a}
\\\exists\ \bmu>0,\ \bw\in U \text{such that}
\mbox{$\bx \in B\bw + \ncone{\bw}{U}$ and
$A^T\bu = \bmu\bw$}. \label{c2 multiplier rule b}
\end{gather}
\end{subequations}
\end{corollary}

\subsubsection{Huber penalty}
\label{HuberSection}

 A popular function in the PLQ class is the Huber penalty~\cite{Huber81}:
\[
  \phi(x) =\!\!  \sup_{w\in[-\kappa,\kappa]^n}\!
  \left[{\ip{x}{w}-\half\tnorm{w}^2}\right] =
      \sum_{i=1}^n \phi_i(x_i); \ \  \phi_i (x_i):=
      \begin{cases}
       \half x_i^2                   &\!\! \mbox{if $|x_i |\le \kappa$,}\\
      \kappa|x_i| - \kappa^2/2 &\!\! \mbox{otherwise.}
    \end{cases}
\]
The Huber function is of form~\eqref{LQ phi}, with $B= I$ and $U = [-\kappa,
\kappa]^{n}$.  In this case, $U^\infty\cap\Nul{B}=\{0\}$ so that the
conditions of Corollary \ref{c2 multiplier cor} hold.

\begin{figure}
\centering
\tikzsetnextfilename{huber_penalty}
\begin{tikzpicture}
  \begin{axis}[
    thick,
    width=.45\textwidth, height=2cm,
    xmin=-2,xmax=2,ymin=0,ymax=1,
    no markers,
    samples=50,
    axis lines*=left, 
    axis lines*=middle, 
    scale only axis,
    xtick={-1,1},
    xticklabels={$-\kappa$,$+\kappa$},
    ytick={0},
    ] 
\addplot[red,domain=-2:-1,densely dashed]{-x-.5};
\addplot[domain=-1:+1]{.5*x^2};
\addplot[red,domain=+1:+2,densely dashed]{x-.5};
\addplot[blue,mark=*,only marks] coordinates {(-1,.5) (1,.5)};
  \end{axis}
\end{tikzpicture}
\tikzsetnextfilename{vapnik_penalty}
\begin{tikzpicture}
  \begin{axis}[
    thick,
    width=.45\textwidth, height=2cm,
    xmin=-2,xmax=2,ymin=0,ymax=1,
    no markers,
    samples=50,
    axis lines*=left, 
    axis lines*=middle, 
    scale only axis,
    xtick={-1,1},
    xticklabels={$-\epsilon$,$+\epsilon$},
    ytick={0},
    ] 
    \addplot[red,domain=-2:-1,densely dashed] {-x-1};
    \addplot[domain=-1:+1] {0};
    \addplot[red,domain=+1:+2,densely dashed] {x-1};
    \addplot[blue,mark=*,only marks] coordinates {(-1,0) (1,0)};
  \end{axis}
\end{tikzpicture}
\caption{Huber (left) and Vapnik (right) penalties}
\label{HuberFigure}
\end{figure}
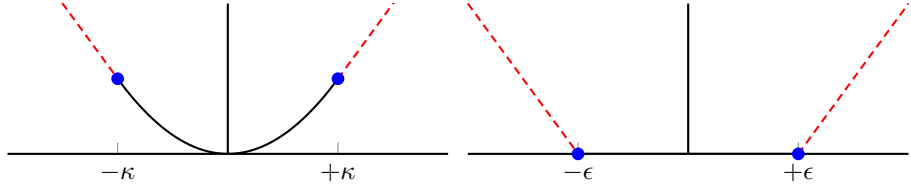

 A graph of the scalar component function $\phi_i$ is shown in Figure~\ref{HuberFigure}. The
Huber penalty is robust to outliers, since it increases linearly
rather than quadratically outside the threshold defined by $\kappa$. For any misfit function $\rho$,
Theorem~\ref{case 2 sd v} can be used to easily compute the
subgradient $\sd v(b,\tau)$ of the value function. If the regularity
condition~\eqref{c2 primal regularity} is satisfied (e.g., if
$\rho$ is finite valued), then Theorem~\ref{case 2 sd v} implies that
\[
\sd v(b,\tau)=\set{\begin{pmatrix}\hfill\bu\\ -\bmu\end{pmatrix} }{
  \begin{array}{c}(\bx,\bu)\mbox{ satisfy \eqref{c2 opt cond} and  }\\
    \bmu=\max\left\{\kappa\|A^T\bu\|_\infty,
      \|A^T\bu\|_2/\sqrt{2\tau}\right\}\end{array}}.
\]
In particular, if $\rho$ is differentiable finite-valued, $\bu =\nabla
\rho(b - A\bx)$ is unique and%
\[
\nabla v(b,\tau) =
\begin{pmatrix}
  \hfill\bu \\ -\bmu
\end{pmatrix}.
\]

\subsection{Affine composition with \RQP  functions} \label{EPLQ}

Next consider the case where $\phi$ takes the form
\begin{equation}\label{eplq}
\phi(x):=\psi(Hx+c),
\quad\mbox{where}\quad
\psi(y):=\sup_{w\in U}[\ip{y}{w}-\half \ip{w}{Bw}],
\end{equation}
$H\in\R^{\nu\times n}$ is injective, $c\in\R^\nu$, $U\subset\R^\nu$ is
nonempty closed and convex with $0\in U$ and $B\in\R^{\nu\times \nu}$
is symmetric and positive semi-definite. We assume that
\[
\exists\, \hx\mbox{ such that }H\hx+c\in\ri{\dom{\psi}},
\]
where
$\dom{\psi}=\cone{U^\circ}+\Ran{B}$ %
(Lemma \ref{c2 phi}).  We show that the function $\phi$ in
\eqref{eplq} is an instance of the quadratic support functions
considered in section~\ref{LinearQuadratic}.  To see this we make the
following definitions:
\[
\begin{array}{c}
\tx=\begin{pmatrix}x\\ s\end{pmatrix},\ 
\ty=\begin{pmatrix}y\\ z\end{pmatrix},\ 
\tw=\begin{pmatrix}v\\ w\end{pmatrix},\ 
\tU\!=\!\{0\}\times U,
\\ \\
\tb\!=\!\begin{pmatrix}b\\ c\end{pmatrix},\  
\tA\!=\!\begin{bmatrix}\phantom-A&0\\ -H&I\end{bmatrix},\  
\tB\!=\!\begin{bmatrix}0&0\\ 0&B\end{bmatrix},\ 
\trho\begin{pmatrix}y\\ z\end{pmatrix}\!=\!\rho(y)+\indicator{z}{\{0\}},
\mbox{ and }
\\ \\
\tphi\begin{pmatrix}x\\ s\end{pmatrix}\!=\!
\displaystyle\sup_{\mbox{\tiny $\begin{pmatrix}v\\ w\end{pmatrix}$}\in\tU}
\left[ \ip{\begin{pmatrix}v\\ w\end{pmatrix}}{\begin{pmatrix}x\\ s\end{pmatrix}}
-\half \ip{\begin{pmatrix}v\\ w\end{pmatrix}}{\tB \begin{pmatrix}v\\ w\end{pmatrix}}\right]
\!=\!\support{x}{\{0\}}+\psi(s).
\end{array}
\]
With these definitions, the two problems \Pbt and
\[\minim \quad \trho(\tb-\tA \tx) \text{subject to} \tphi(\tx)\le\tau
\]
are equivalent.  In addition, we have the relationships
\begin{equation*}
\begin{array}{c}
\trho^*\begin{pmatrix}u\\ r\end{pmatrix}=\rho^*(u)+\support{r}{\{0\}},
\quad
\tphi^*\begin{pmatrix}v\\ w\end{pmatrix}=\indicator{v}{\{0\}}+\psi^*(w),\\ \\
\gauge{\begin{pmatrix}v\\ w\end{pmatrix}}{\tU}=\indicator{v}{\{0\}}+\gauge{w}{U},
\text{and}
\norm{\begin{pmatrix}v\\ w\end{pmatrix}}_{\tB}=\support{v}{\{0\}}+\norm{w}_B.
\end{array}
\end{equation*}
Moreover, the reduced dual \ref{GeneralValueDualr} becomes
\begin{equation}
\sup_{H^Tr=A^Tu}\left[\ip{b}{u}+\ip{c}{r}-\rho^*(u)-\support{r}{\lev{\psi}{\tau}}\right].
\label{c3 reduced dual}
\end{equation}
Using standard methods of convex analysis, we obtain the following
result as a direct consequence of Theorem~\ref{case 2 sd v}
and~\cite[Corollary~10.11]{RTRW}.

\begin{theorem}\label{c3 sd v}
Let
$\phi$ be given by \eqref{eplq}, and
consider the following two conditions:
\begin{equation}\label{c3 primal regularity}
\exists\ \hx\mbox{ such that }H\hx+c\in\ri{\dom{\psi}},\ \psi(H\hx +c)<\tau,\mbox{ and }b-A\hx\in\ri{\dom{\rho}}
\end{equation}
and
\begin{equation}\label{c3 dual regularity}
\exists\ \hu\in\ri{\dom{\rho^*}}\mbox{ and }\hr\in\ri{\cone{U}}\mbox{  such that  }
\begin{pmatrix}\hu\\ \hr\end{pmatrix}\in\Nul{\begin{bmatrix}\hfill A\\ -H\end{bmatrix}^T}.
\end{equation}

\medskip

\noindent 1.
If the triple $(\bx,\bu,\br)$ satisfies
\begin{equation}\label{c3 opt cond}
\bx\in\lev{\phi}{\tau},\ 
\bu\in\sd\rho(b-A\bx),\ \br\in\ncone{H\bx+c}{\lev{\psi}{\tau}}
\mbox{ and }  A^T\bu=H^T\br\ ,
\end{equation}
then $\bx$ solves  $\Pbt$ and $(\bu,\br)$ solves \eqref{c3 reduced dual}.

\medskip

\noindent 2.  If $\bx$ solves $\Pbt$ and \eqref{c3 primal regularity}
holds, there exists $(\bu,\br)$ such that \eqref{c3 opt cond} holds.

\medskip

\noindent 3. If $(\bu,\br)$ solves \eqref{c3 reduced dual} and
\eqref{c3 dual regularity} holds, there exists $\bx$ such that
\eqref{c3 opt cond} holds.

\medskip

\noindent 4.
If
either
\[
H^{-1}[\cone{U}^\circ]\cap[-A^{-1}\hzn{\rho}]=\{0\}\mbox{ and
\eqref{c3 primal regularity} holds, }
\]
or
\[
\begin{pmatrix}b\\ c\end{pmatrix}\in
\intr{\dom{\rho}\times\lev{\psi}{\tau}+\Ran{\begin{bmatrix}\hfill A\\ -H\end{bmatrix}}}
\mbox{ and \eqref{c3 dual regularity} holds,}
\]
then $\sd v(b,\tau)\ne\emptyset$ and is given by
\begin{align*}%
\sd v(b,c,\tau)&=\left\{{\begin{pmatrix}\hfill\bu\\\hfill\br\\-\bmu\end{pmatrix}}\left\vert
\begin{array}{c}\exists\; \bx\in\Rn\mbox{ s.t. } (\bx,\bu,\br)\mbox{ satisfy \eqref{c3 opt cond} and  }\\
\bmu=\max\left\{\gauge{\br}{U}, \|\br\|_B/\sqrt{2\tau}\right\}\end{array}\right.\right\}
\\
&=\left\{\begin{pmatrix}\hfill\bu\\\hfill\br\\-\bmu\end{pmatrix}\left\vert
\begin{array}{c}
\exists\; \bx\in\Rn\mbox{ s.t. } c+H\bx\in\lev{\psi}{\tau},\\
\bu\in\sd \rho(b-A\bx),\ \br\in\bmu^+\sd \psi(c+H\bx), \ \bmu\ge 0,\\
\bmu(\psi(c+H\bx)-\tau)=0,\mbox{ and }A^T\bu=H^T\br
\end{array}
\right.\right\}.
\end{align*}

\end{theorem}

\begin{corollary}\label{c3 multiplier cor}
  Consider the problem $\Pbt$ with $\phi$ given by \eqref{eplq}.  Then
  $(\bx,\bu,\br)$ satisfies \eqref{c3 opt cond} if and only if
  \begin{gather*}
    H\bx + c \in \lev{\psi}{\tau}, \quad
    \bu \in \sd\rho(b-A\bx), \quad
    A^T \bu = H^T \br,
    \\\text{and either}
    \br\in\ncone{H\bx+ c}{\dom\psi}, \text{or}
    \\\exists\ \bmu\ge0,\ \bw\in U \text{such that}
    \mbox{$Hx + c\in B\bw + \ncone{\bw}{U}$ and $\br=\bmu\bw$.}
  \end{gather*}
\end{corollary}

\subsubsection{Vapnik penalty}

The Vapnik penalty 
\[
\rho(r)
=
\sup_{u \in [0,1]^{2n}}
\left\{
\left\langle \begin{bmatrix}\phantom-r - \epsilon \\
-r - \epsilon \end{bmatrix} , u \right\rangle
\right\}
=
(r - \epsilon)_+ + (-r-\epsilon)_+
\]
is an important example in the PLQ class which 
is most easily represented as the composition of an affine transformation
with a PLQ function.
The scalar
version is shown in the right panel of Figure~\ref{HuberFigure}.
In this case,
\[
H = \begin{bmatrix}\phantom-I \\ -I\end{bmatrix},\quad
c =  - \begin{bmatrix} \epsilon \B{1}\\ \epsilon \B{1}\end{bmatrix},\quad
B= \B{0} \in \bR^{2n \times 2n}, \text{and} U = [0, 1]^{2n}.
\]
In order to satisfy~\eqref{c3 opt cond}, we need to find a triple
$(\bx, \bu, \bw)$ with $\bw = [\bw_1 \ \bw_2]^T\in [0,1]^{2n}$ so that
$\bu \in \partial \rho(b - A\bx)$ and $A^T\bu = H^T\bw = \bw_1 -
\bw_2$.  We claim that either $\bw_1(i)= 0$ or $\bw_2(i) = 0$ for all
$i$.  To see this, observe that
$\bw\in\ncone{H\bx+c}{\lev{\psi}{\tau}}$, so
\[
\left\langle \bw, y - \begin{bmatrix} \phantom-\bx - \epsilon \\ -\bx - \epsilon\end{bmatrix}\right\rangle \leq 0
\]
whenever $\psi(y) \leq \tau$. Taking $y$ first with $-\epsilon$ as the only non-zero in the
$i$th coordinate, and then with $-\epsilon$ in the only nonzero in the $(n+i)$th coordinate, we get
\[
\bw_1(i) (-\bx(i)) \leq 0 \text{and} \bw_2(i) (\bx(i)) \leq 0.
\]
If $x(i) < 0$, from the first equation we get $\bw_1(i) = 0$, while if
$x(i) > 0$, we get $\bw_2(i) = 0$ from the second equation.  If $x(i)
= 0$, then taking $y = 0$ gives
\[
\bw_1(i) \epsilon \leq 0 \text{and} \bw_2(i) \epsilon \leq 0,
\]
so $\bw_1(i) = \bw_2(i) = 0$.
Since $A^T\bu = \bw_1 - \bw_2$, and $\bw_1(i)$ or $\bw_2(i)$ is  $0$ for each $i$, we get $\mu = \gauge{\bw}{[0,1]^{2n}} = \|A^T\bu\|_\infty$.
Hence, the subdifferential $\partial v$ is computed in precisely the same way for the Vapnik regularization as for the $1$-norm.

\section{Numerical example: robust nonnegative BPDN}
\label{Numerics}

In this example, we recover a nonnegative undersampled sparse signal
from a set of very noisy measurements using several formulations
of~\ref{ConstrainedClass}.  We compare the performance of three
different penalty functions $\rho$: least-squares, Huber (see
section~\ref{HuberSection}), and a nonconvex penalty arising from the
Student's t distribution (see, e.g.,~\cite{AravkinFHV:2012,SYSID2012tks}).  The
regularizing function $\phi$ in all of the examples is the sum of the
1-norm and the indicator of the positive orthant (see
section~\ref{NormIndicatorSection}).

The formulations using Huber and Student's t misfits are 
robust alternatives to the nonnegative basis pursuit problem~\cite{DON2005Ta}. 
The Huber misfit agrees with the quadratic penalty for small residuals, but
is relatively insensitive to larger residuals. The Student's t misfit 
is the negative likelihood of the Student's t distribution, 
\begin{equation}
\label{student}
\rho_s(x) = \log(1 + x^2/\nu),
\end{equation}
where $\nu$ is the degrees of freedom parameter. 

For each penalty $\rho$, our aim is to solve the problem
\[
 \minim_{x\ge0} \quad \norm{x}_1 \quad\st\quad\rho(b - Ax)\le\sigma,
\]
via a series of approximate solutions of~\ref{ConstrainedClass}. The
1-norm regularizer on $x$ encourages a sparse solution.  In
particular, we solve the nonlinear equation~\eqref{eq:v-tau-sig},
where $v$ is the value function of~\ref{ConstrainedClass}. This is the
approach used by the SPGL1 software
package~\cite{BergFriedlander:2011}; the underlying theory, however,
does not cover the Huber function.  Also, $\phi$ is not everywhere
finite valued, which violates~\cite[Assumption
3.1]{BergFriedlander:2011}.  Finally, the Student's t
misfit~\eqref{student} is nonconvex; however, the inverse function
relationship (cf.~Theorem~\ref{IVT}) still holds, so we can achieve
our goal, provided we can solve the root-finding problem.

Formula~\eqref{HuberTauDerivative} computes the derivative of the
value function associated with \Pbt for any convex differentiable
$\rho$. The derivative requires $\nabla \rho$, evaluated at the
optimal residual associated with \Pbt.  For the Huber case, this is
given by
\[
(\nabla \rho(b - A\bar x))_i
=
\sign(b_i - A_i\bar x)\cdot\min(| b_i - A_i\bar x|, \kappa).
\]
The Student's t misfit is also smooth, but nonconvex. Therefore, the
formula~\eqref{HuberTauDerivative} may still be applied---with the
caveat that there is no guarantee of success. However, in all of the
numerical experiments, we are able to find the root
of~\eqref{eq:v-tau-sig}.

We consider a common compressive sensing example: we want to recover
a 20-sparse vector in $\R_+^{512}$ from 120 measurements. We use a Gaussian measurement
matrix $A \in \R^{100 \times 1024}$, where each entry is sampled from
the distribution $N(0, 1/10)$. We generate measurements to test the BPDN formulation
according to
\[
b = Ax + w + \zeta,
\]
where $w \sim N(0, 0.005^2)$ is small Gaussian error, 
while $\zeta$ contains five randomly placed large outliers  
sampled from $N(0, 4)$.  For each penalty $\rho$, the 
$\sigma$ parameter is the true measure of the error in that penalty, 
i.e., $\sigma_\rho = \rho(\zeta)$. This allows a fair comparison between 
the penalties.

\begin{figure}
  \centering
  { \includegraphics[width=.48\textwidth]{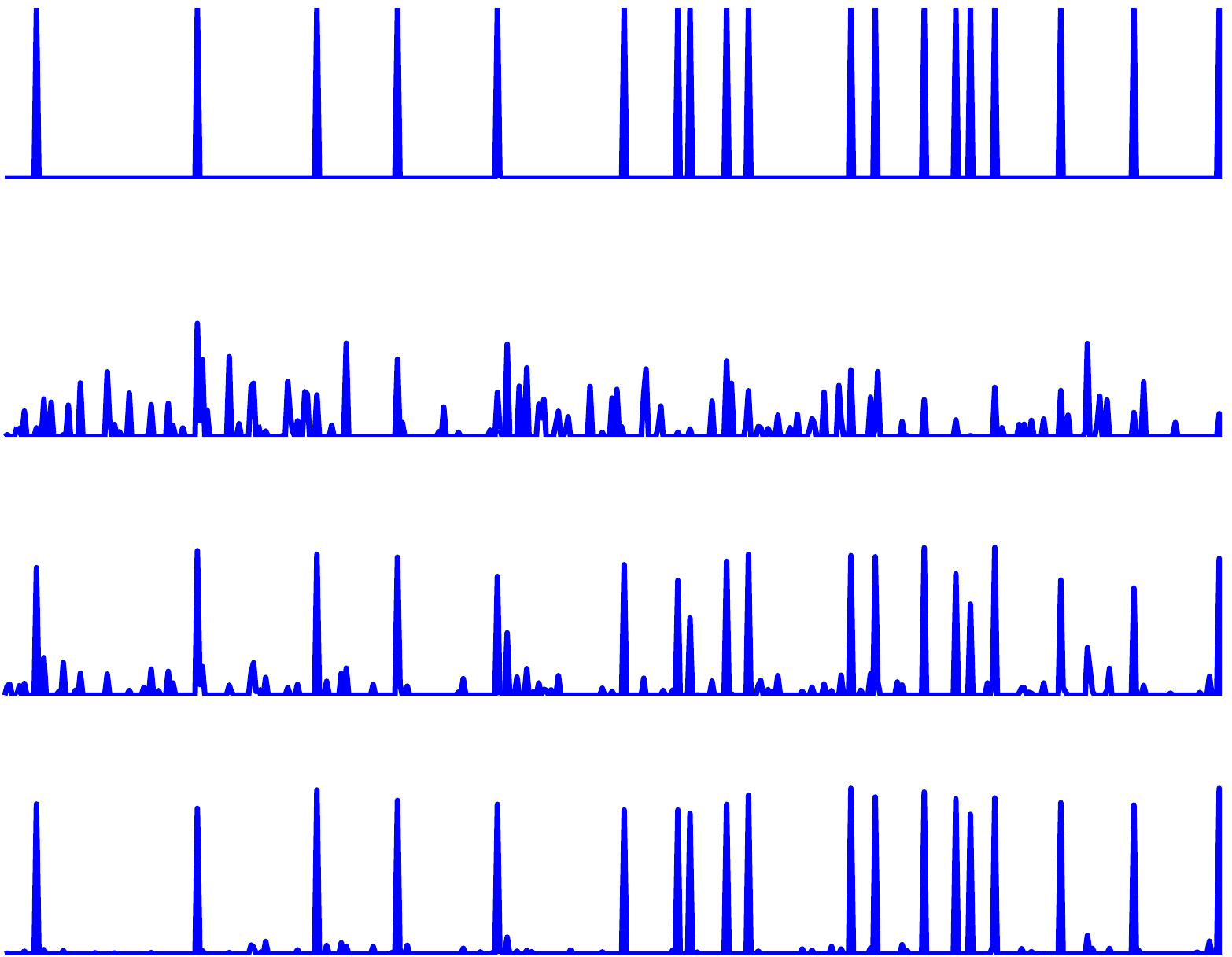}}
\raisebox{-2ex}{
  {\includegraphics[width=.48\textwidth]{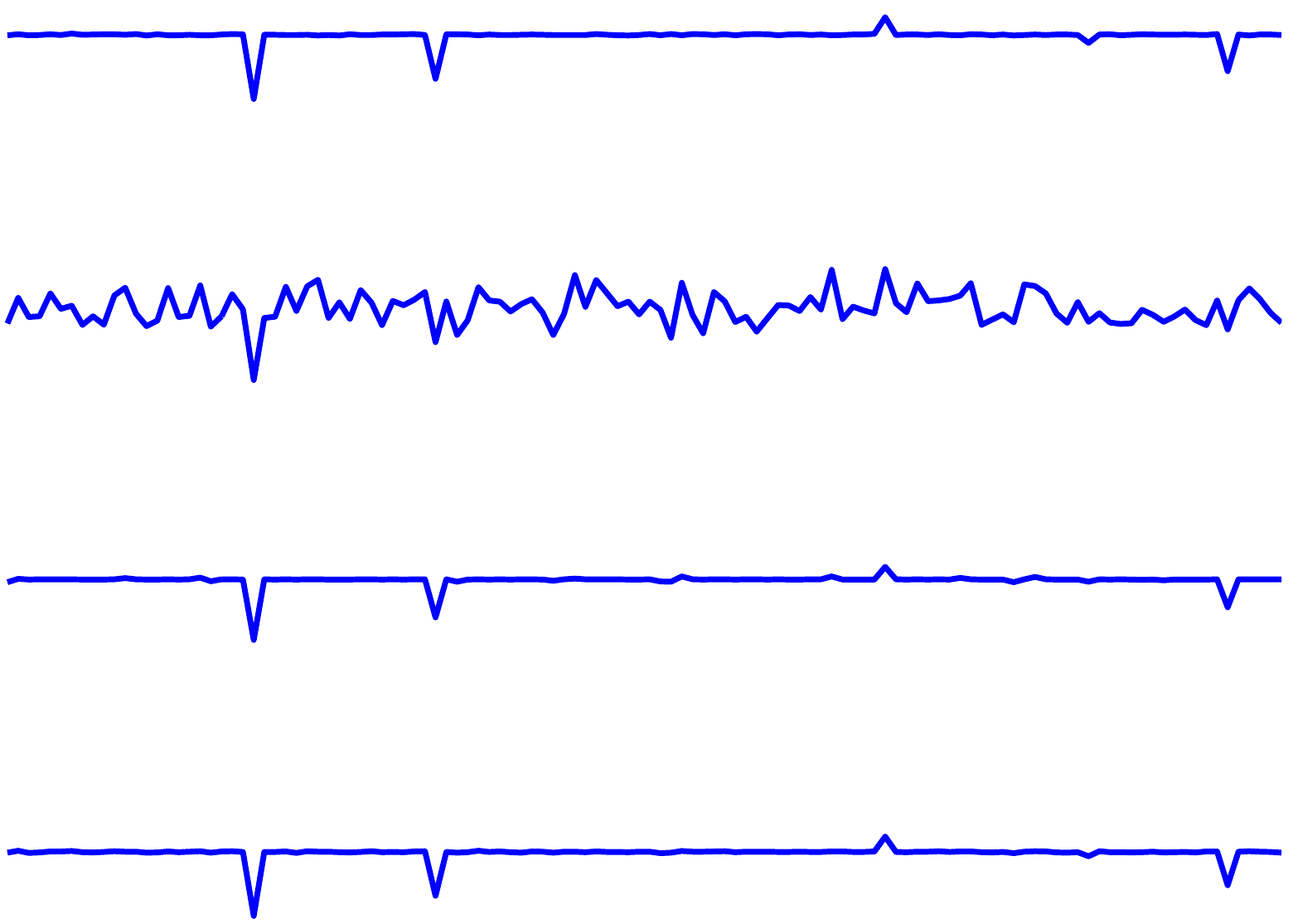}}
  }
  \caption{\label{Results} Left, top to bottom: True signal, and
    reconstructions via least-squares, Huber, and Student's t.  Right,
    top to bottom:  true errors, and least-squares,
    Huber, and
    Student's t residuals.}
\end{figure}

We expect the Huber function to out-perform the least squares penalty
by budgeting the error level $\sigma$ to allow a few large outliers,
which will never happen with the quadratic. We expect the Student's t
penalty to work even better, because it is non convex, and the grows
sublinearly as outliers increase.  The results in Figure~\ref{Results}
demonstrate that this is indeed the case.  In many instances the Huber
function is able to do just as well as the Student's t; however, often
the Student's t does better (and never worse). Both robust penalties
always do better than the least squares fit.  The code is implemented
in and extended version of SPGL1, and can be downloaded from
\verb{https://github.com/saravkin/spgl1{.  The particular experiment
    presented here can be found in \verb{tests/spgl1TestNN.m{.

\section{Appendix: Proofs of results} \label{sec:proofs-results}

\subsection*{Proof of Theorem~\ref{IVT}}

Let $\tau\in\cS_{1,2}$ and set $\sigma_\tau=v_1(\tau)$.  By assumption,
$\argmin \cP_{1,2}(\tau)\ne\emptyset$.  Let $x_\tau\in\argmin
\cP_{1,2}(\tau)$, so that $\psi_1(x_\tau)=\sigma_\tau$ and
$\psi_2(x_\tau)=\tau$. In particular, $x_\tau$ is feasible for
$\cP_{2,1}(\sigma_\tau)$. Let $\hx$ be any other feasible point for
$\cP_{2,1}(\sigma_\tau)$ so that $\psi_1(\hx)\le\sigma_\tau=
v_1(\tau)=\psi_1(x_\tau)$.  If $\psi_1(\hx)<\sigma_\tau=v_1(\tau)$, then
$\psi_2(\hx)>\tau$ since otherwise we contradict the definition of
$v_1(\tau)$. If $\psi_1(\hx)=\sigma_\tau$, then we claim that
$\psi_2(\hx)\ge\tau$. Indeed, if $\psi_2(\hx)<\tau$, then
$\hx\in\argmin \cP_{1,2}(\tau)$ but $\psi_2(\hx)<\tau$, which
contradicts the fact that $\tau\in\cS_{1,2}$. Hence, every feasible
point for $\cP_{2,1}(\sigma_\tau)$ has $\psi_2(\hx)\ge\tau$ with
equality only if $\psi_1(\hx)=\sigma_\tau$. But $x_\tau$ is feasible for
$\cP_{2,1}(\sigma_\tau)$ with $\psi_2(x_\tau)=\tau$. Therefore,
$x_\tau\in \argmin\cP_{2,1}(\sigma_\tau)\subset\set{x\in
  X}{\psi_1(x)=\sigma_\tau}$. Consequently, $v_2(v_1(\tau))=\tau$ and
\begin{equation}\label{inverse inclusion}
\emptyset\ne\argmin\cP_{1,2}(\tau)\subset\argmin\cP_{2,1}(\sigma_\tau)\subset\set{x\in X}{\psi_1(x)=\sigma_\tau}.
\end{equation}

We now show that
$\argmin\cP_{2,1}(\sigma_\tau)\subset\argmin\cP_{1,2}(\tau)$. Let
$\hx\in \argmin\cP_{2,1}(\sigma_\tau)$.  In particular, $\hx$ is
feasible for $\cP_{2,1}(\sigma_\tau)$, so, by what we have already
shown, $\psi_2(\hx)\ge\tau$ with equality only if
$\psi_1(\hx)=\sigma_\tau$. But, by our choice of $\hx$,
$\psi_2(\hx)=v_2(v_1(\tau))=\tau$, so $\psi_1(\hx)=\sigma_\tau$, i.e.,
$\hx\in \argmin\cP_{1,2}(\tau)$.

It remains to establish the final statement of the theorem.  By
\eqref{inverse inclusion}, we already have that
$\set{v_1(\tau)}{\tau\in\cS_{1,2}}\subset\cS_{2,1}$, so we need only
establish the reverse inclusion. For this, let $\sigma\in\cS_{2,1}$ and
set $\tau_\sigma=v_2(\sigma)$. By interchanging the indices and applying
the first part of the theorem, we have from \eqref{inverse inclusion}
that
\begin{equation*}
\emptyset\ne\argmin\cP_{2,1}(\sigma)\subset\argmin\cP_{1,2}(\tau_\sigma)\subset\set{x\in X}{\psi_2(x)=\tau_\sigma}.
\end{equation*}
That is, $\tau_\sigma\in \cS_{1,2}$ and, by (a), $\sigma=v_1(v_2(\sigma))=v_1(\tau_\sigma)$.

\subsection*{Proof of the inverse linear image
  (section~\ref{sec:calculus})}
 \label{app:inverse linear image}

For $\lam>0$, observe that
\begin{eqnarray}
h^{\pi}(w,\lam)&=&\lam \inf_{A x=\lam^{-1}w}p(x)\nonumber\\
&=&\lam \inf_{A s=w}p(\lam^{-1}s)\qquad(s:=\lam x)\nonumber\\
&=&\inf_{A s=w}p^\pi(s,\lam)\label{iac2}\\
&=&\inf\set{p^\pi(s,\zeta) }{
\hat A{ \left(\begin{array}{c}s\\ \zeta\end{array}\right)
=
\left(\begin{array}{c}w\\ \lam\end{array}\right)}},\label{iac3}
\end{eqnarray}
where
\[
\hat A=\left[\begin{array}{cc}A&0\\ 0&1\end{array}\right].
\]
Again by \cite[Theorem 9.2]{RTR} in conjunction with
\cite[Corollary 16.2.1]{RTR}, the function in \eqref{iac3} is closed if
$(\hat A^T)^{-1}\dom{(p^\pi)^*}\ne\emptyset$. Since, by~\cite[Corollary 13.5.1]{RTR},
$\dom{(p^\pi)^*}=\set{(u,\eta) }{ p^*(u)\le -\eta}$, we have
 \[
 (\hat A^T)^{-1}\dom{(p^\pi)^*}\ne\emptyset
 \text{if and only if}
(A^T)^{-1}\dom{p^*}\ne\emptyset.
\]
Hence, by assumption, the function in \eqref{iac3}
is closed proper convex and equals $h^{\pi}(w,\lam)$ on the relative interior
of its domain. Since $h^{\pi}(w,\lam)$ is closed, \eqref{iac2} implies that these functions must coincide.

\subsection*{Proof of Lemma~\ref{lem:conjugate of epi indicator}}

We first prove~\eqref{epi support and perspective}.  The conjugate
of $\indicator{(x,\tau)}{\epi{h}}$ is obtained as follows:
\begin{align*}
\support{(y,\mu)}{\epi{h}}
&=
\sup_{\tau,x}\left[ \ip{y}{x}+\mu\tau-\indicator{(x,\tau)}{\epi{h}}\right]\\
&=
\sup_{\tau,x\in\dom{h}}\left[ \ip{y}{x}+\mu\tau-\indicator{h(x)-\tau}{\bR_-}\right]\\
&=
\sup_{\omega,x\in\dom{h}}
\!\!\big[ \ip{y}{x}+\mu(h(x)-\omega)-\indicator{\omega}{\bR_-}\big]\quad\quad(\omega:=h(x)-\tau)\\
&=
\sup_{x\in\dom{h}}
\!\big[ \ip{y}{x}+\mu h(x)+\sup_\omega\left[-\mu\omega-\indicator{\omega}{\bR_-}\right]\big]\\
&=
\sup_{x\in\dom{h}}
\left[ \ip{y}{x}+\mu h(x)+\indicator{\mu}{\R_-}\right].
\end{align*}
For $\mu<0$, we obtain
\[
\support{(y,\mu)}{\epi{h}}
   = -\mu\sup_{x}\left[ \ip{-\mu^{-1}y}{x}-h(x)\right] 
= -\mu h^*(-\mu^{-1}y).
\]
Since $h^*$ is necessarily a closed proper convex function, we obtain
the result.

To see \eqref{level support}, first note that the function
\[
q(y):=\inf_{\mu> 0}[\tau\mu+\mu h^*(y/\mu)]=\inf_{\mu\ge 0}[\tau\mu+(h^*)^\pi(y,\mu)]
\]
is the positively homogeneous function generated by the function
$y\mapsto \tau+h(y)$~\cite[page 35]{RTR}, and so is convex in $y$.
Next observe that the conjugate of $q$ is given by
\begin{eqnarray*}
q^*(x)&=&\sup_{y}\left[\ip{x}{y}-\inf_{\mu\ge 0}[\tau\mu+(h^*)^\pi(y,\mu)]\right]
\\ &=&
\sup_{y,\mu\ge 0}\left[\ip{x}{y}+\tau(-\mu)-(h^*)^\pi(y,\mu)\right]
\\ &=&
\sup_{(y,\mu)}\left[\ip{x}{y}+\tau\mu-(h^*)^\pi(y,-\mu)\right]\quad\mbox{(exchange $-\mu$ for $\mu$)}
\\ &=&
\sup_{(y,\mu)}\left[\ip{x}{y}+\tau\mu-\support{(y,\mu)}{\epi{h}}\right]\quad
\mbox{(by~\eqref{epi support and perspective})}
\\ &=&
\indicator{(x,\tau)}{\epi{h}}
=%
\indicator{x}{\lev{h}{\tau}}.
\end{eqnarray*}
The result now follows from the Bi-Conjugate Theorem~\cite[Theorem 12.2]{RTR}.

\subsection*{Proof of Theorem~\ref{reduced dual theorem}}

Combining \ref{dual1} with \eqref{level support} and \eqref{f conj} gives
\begin{align}
\hv(b,\tau)&:=\sup_{u,\mu}\left[\ip{b}{u}+\tau\mu- (\phi^*)^{\pi}(A^Tu,-\mu)-\rho^*(u)\right]
\label{GeneralValueDual}\\ &=
\sup_{u}\left[\ip{b}{u}-\rho^*(u)-\inf_{\mu\le 0}\left[\tau(-\mu)+(\phi^*)^{\pi}(A^Tu,-\mu)\right]\right]
\nonumber\\ &=
\sup_{u}\left[\ip{b}{u}-\rho^*(u)-\inf_{\mu\ge 0}\left[\tau(\mu)+(\phi^*)^{\pi}(A^Tu,\mu)\right]\right]
\nonumber\\ &=
\sup_{u}\left[\ip{b}{u}-\rho^*(u)-\support{A^Tu}{\lev{\phi}{\tau}}\right],
\nonumber
\end{align}
where the final equality follows from \eqref{level support}.
The equivalence \eqref{g conj} follows from the definition of the conjugate and 
the equivalence \eqref{g conj2} follows from \cite[Theorems 16.3 and 16.4]{RTR} 
which tell us that
\begin{eqnarray*}
g_\tau^*(b)&=&\cl{\rho\infc \left[\support{A^T\cdot}{\lev{\phi}{\tau}}\right]^*}(b)
\\ &=&\cl{\inf_{w_1+w_2=\cdot} \left[ \rho(w_1)+\inf_{Ax=w_2}\indicator{x}{\lev{\phi}{\tau}}
\right]}(b)
\\ &=&\cl{\inf_{\phi(x)\le\tau}\rho(\cdot-Ax)}(b)
\\ &=&\cl{v(\cdot,\tau)}(b)\, .
\end{eqnarray*}
The uniqueness of $u$ when $\rho$ is differentiable follows from the
essential strict convexity of $\rho^*$~\cite[Theorem~26.3]{RTR}.

\subsection*{Proof of Lemma~\ref{mu existence}}

\paragraph{Part 1} The inequality follows immediately
from~\eqref{level support}. But it is also easily derived from the
observation that if $\mu>0$ and $x \in \lev{\phi}{\tau}$, then
\[
\begin{aligned}
\tau \mu + \mu \phi^*(s/\mu) & \geq \tau \mu + \mu[\ip{x}{s/\mu} -\phi(x)] 
&\mbox{[Fenchel-Young inequality]}\\
& \geq \phi(x) \mu + \ip{x}{s} - \mu \phi(x) \\
&= \ip{x}{s}.
\end{aligned}
\]
Taking the sup over $x \in \lev{\phi}{\tau}$ gives the result.
\smallskip

\paragraph{Part 2} The proof uses the following three key facts:
\begin{enumerate}
\item[(i)]
By \cite[Theorems 23.5 and 23.7]{RTR},
for any non-empty closed convex set $U$ and $\bu\in U$, 
\begin{equation}\label{support ncone}
\bv \in \ncone{\bu}{U}\ \iff\ 
\bu \in \partial\support{\bv}{U}=\argmax_{u\in U}
\ip{\bv}{u}.
\end{equation}
\item[(ii)]
The Fenchel-Young inequality tells us that
\begin{equation}\label{tau plus phi conj}
\tau+\phi^*(\bs)\ge \phi(\bx)+\phi^*(\bs)\ge\ip{\bs}{\bx}.
\end{equation}
\item[(iii)] (\!\cite[Lemma 26.17]{BauCom:2011} or \cite[Corollary 2.9.5]{Zal:2002})
Let $\map{g}{\eR}$ be a convex function and $\tau\in\R$ be such that
$\tau>\inf g$. Then for every $x\in\lev{g}{\tau}$
\begin{equation}\label{lev ncone}
\ncone{x}{\lev{g}{\tau}}=
\begin{cases}
\ncone{x}{\dom{g}}\cup\cone{\sd g(x)},&\mbox{ if $g(x)=\tau$,}\\
\ncone{x}{\dom{g}},&\mbox{ if $g(x)<\tau$.}
\end{cases}
\end{equation}
\end{enumerate}

\noindent
We divide the proof into two parts: 
(A) if $S_1\ne \emptyset$, show $S_1\subset S_2$, and
(B) if $S_2\ne \emptyset$, show $S_2\subset S_1$ and equality holds in \eqref{LevelSupportInf}.
Combined, these implications establish Part 2 of the lemma.

\noindent
(A) Let $\bmu\in S_1$. We show that $\bmu\in S_2$.
First suppose $\phi(\bx)<\tau$. 
By \eqref{lev ncone},
$\ncone{\bx}{\lev{\phi}{\tau}}=\ncone{\bx}{\dom{\phi}}$.
Hence,  
by \eqref{support ncone}, $\bs\in \ncone{\bx}{\dom{\phi}}$. Therefore, if
$\bmu=0$, we have $\bmu\in S_2$. 

On the other hand, if $\bmu>0$, 
by \eqref{support ncone} and the fact that $\ncone{\bx}{\lev{\phi}{\tau}}=\ncone{\bx}{\dom{\phi}}$, we have
\[
\begin{aligned}
\ip{\bs}{\bx}&=\support{\bs}{\dom{\phi}}\\
&=(\phi^*)^\infty (\bs)&\mbox{\cite[Theorem 13.3]{RTR}}\\
&=\tau 0 +(\phi^*)^\pi(\bs,0)\\
&\ge \tau \bmu +(\phi^*)^\pi(\bs,\bmu)\\
&> \bmu\phi(\bx)+\bmu \phi^*(\bs/\bmu)\\
&\ge \bmu\ip{\bs/\bmu}{\bx}&\mbox{[Fenchel-Young inequality]}\\
&=\ip{\bs}{\bx}.
\end{aligned}
\]
Since this cannot occur, it must be the case that $\bmu=0$ and $\bmu\in S_2$.

Now suppose that $\phi(\bx)=\tau$ and $\bs=0$.
Then, for $\mu> 0$, 
$p_\tau(\bs,\mu)=(\tau+\phi^*(0))\mu\ge 0$ by \eqref{tau plus phi conj}, and,
for $\mu =0$,  $p_\tau(\bs,\mu)=(\phi^*)^\infty (0)=0$.
Therefore,
$0=\inf_{0\le\mu}p_\tau(\bs,\mu)$ with $\mu =0\in S_1$. But, in this
case, it is also clear $0\in S_2\ne\emptyset$, since $\bs = 0 \in \ncone{\bx}{\dom{\phi}}$. Thus, if $\bmu =0$ we have $\bmu\in S_2$.
If $\bmu>0$, then $0=\tau+\phi^*(0)$ since $0=p_\tau(0,\bmu)=(\tau+\phi^*(0))\bmu$. 
But then, by \eqref{tau plus phi conj}, $\phi(\bx)+\phi^*(0)=\ip{\bs}{\bx}=0$ so that $\bs=0\in \sd\phi(\bx)$.
However, $\phi(\bx)=\tau>\inf \phi$, so $0\notin \sd\phi(\bx)$
\cite[Theorem 23.5(b)]{RTR}. This contradiction implies that if
$\bs =0$, then we must also have $\bmu=0$, and, in particular, we have $S_1\subset S_2$.

Finally, suppose that $\phi(\bx)=\tau$ and $\bs\ne 0$. 
Then, by \eqref{lev ncone}, 
\begin{equation*} 
\mbox{either \quad (a) $\bs\in\cone{\sd \phi(\bx)}$ \quad or \quad (b) $\bs\in\ncone{\bx}{\dom{\phi}}$.} 
\end{equation*}
Let us first suppose that $\bs\notin\ncone{\bx}{\dom{\phi}}$ so, in particular, $\bs\in\cone{\sd \phi(\bx)}$.
As an immediate consequence, we have that $S_2\ne\emptyset$ and the only values of $\mu$ for
which $\bs\in \mu^+\sd\phi(\bx)$ have $\mu>0$ since $\bs\notin\ncone{\bx}{\dom{\phi}}$. 
Let 
$0<\hat\mu\in S_2$. If $\bmu=0$, then
\[
\begin{aligned}
\support{\bs}{\dom{\phi}}&=(\phi^*)^\infty(\bs)\\
&=\inf_{0\le\mu}p_\tau(\bs,\mu)\\
&\le \tau\hat\mu+\hat\mu\phi^*(\bs/\hat\mu)\\
&=\hat\mu\phi(\bx)+\hat\mu[\ip{\bs/\hat\mu}{\bx}-\phi(\bx)]&
\left[\begin{array}{c}\bs/\hat\mu\in\sd \phi(\bx)\mbox{ and}\\ \mbox{\cite[Theorem 23.5(d)]{RTR}}
\end{array}\right]
\\
&=\ip{\bs}{\bx}\\
&\le \support{\bs}{\dom{\phi}}
\end{aligned}
\]
so that $\ip{\bs}{\bx}=\support{\bs}{\dom{\phi}}$, or equivalently, $\bs\in\ncone{\bx}{\dom{\phi}}$
contradicting the choice of $\bs$. Hence, it must be the case that $\bmu>0$.
Again let $0<\hat\mu\in S_2$. Then, by Part 1,
\[
\begin{aligned}
\support{\bs}{\lev{\phi}{\tau}}&\le p_\tau(\bs,\bmu)\\
&=\inf_{0\le\mu}p_\tau(\bs,\mu)\\
&\le \tau\hat\mu+\hat\mu\phi^*(\bs/\hat\mu)\\
&=\hat\mu\phi(\bx)+\hat\mu[\ip{\bs/\hat\mu}{\bx}-\phi(\bx)]&\left[\begin{array}{c}\bs/\hat\mu\in\sd \phi(\bx)\mbox{ and}\\ \mbox{\cite[Theorem 23.5(d)]{RTR}}
\end{array}\right]
\\
&=\ip{\bs}{\bx}\\
&\le \support{\bs}{\lev{\phi}{\tau}}
\end{aligned}
\]
so that $\ip{\bs}{\bx}=\bmu[\phi(\bx)+\phi^*(\bs/\bmu)]$, or equivalently, $\bs\in\bmu\sd\phi(\bx)$.
Hence, $\bmu\in S_2$.

Finally, consider the case where $0\ne \bs\in\ncone{\bx}{\dom{\phi}}$. Then
\[
\begin{aligned}
\inf_{\mu\ge 0} p_\tau(\bs,\mu)&\le p_\tau(\bs,0)
\\
&=(\phi^*)^\infty(\bs)
\\
&=\support{\bs}{\dom{\phi}}&\quad&\mbox{\cite[Theorem 13.3]{RTR}}
\\
&=\ip{\bs}{\bx}&\quad&\mbox{[by \eqref{support ncone}]}
\\
&=\support{\bs}{\lev{\phi}{\tau}}&\quad&\mbox{[again by \eqref{support ncone}]}
\\
&\le \inf_{\mu\ge 0} p_\tau(\bs,\mu)\; ,&\quad&\mbox{[Part 1]}
\end{aligned}
\]
so $0\in S_1$ and $0\in S_2$. If $\bmu>0$, then this string of equivalences also
implies that $\ip{\bs}{\bx}=p_\tau(\bs,\bmu)=\bmu[\phi(\bx)+\phi^*(\bs/\bmu)]$,
or equivalently, $\bs\in\bmu\sd\phi(\bx)$ so that $\bmu\in S_2$.
Putting this all together, we get that $S_1\subset S_2$.

\noindent
(B) Let $\bmu\in S_2$. If $\bmu =0$, then
\[
\begin{aligned}
p_\tau(\bs,0) %
&=(\phi^*)^\infty(\bs) 
\\
&=\support{\bs}{\dom \phi} &\mbox{\cite[Theorem 13.3]{RTR}}
\\
&=\ip{\bs}{\bx}
\\
&\le\support{\bs}{\lev{\phi}{\tau}}
\\
&\le\inf_{\mu\ge 0}\, p_\tau(\bs,\mu).&\mbox{[Part 1]}
\end{aligned}
\]
Therefore, $\bmu=0\in S_1$ and equality holds in \eqref{LevelSupportInf}.

On the other hand, if $\bmu>0$, then
$\bs/\bmu\in\partial\phi(\bx)$, and so
\[ 
\begin{aligned}
\tau\bar\mu+(\phi^*)^\pi(\bs,\bar\mu)
&=
\bar\mu[\phi(\bx)+\phi^*(\bs/\bar\mu)]\\
&=\bar\mu\ip{\bx}{\bs/\bar\mu}&\left[\begin{array}{c}\bs/\hat\mu\in\sd \phi(\bx)\mbox{ and}\\ \mbox{\cite[Theorem 23.5(d)]{RTR}}
\end{array}\right]
\\
&=\ip{\bx}{\bs}\\
&\le\support{\bs}{\lev{\phi}{\tau}}\\
&\le\inf_{\mu\ge 0}\left[\tau\mu+(\phi^*)^{\pi}(\bs,\mu)\right].&\mbox{[Part 1]}
\end{aligned}
\]
Hence, $\bmu\in S_1$ and equality holds in \eqref{LevelSupportInf}.

\subsection*{Proof of Lemma~\ref{coercivity 1}}
\label{sec:proof-lemma-refc}

\paragraph{Part 1} 
The primal coercivity equivalence follows from \cite[Theorems 8.4 and 8.7]{RTR}
since $\hzn{f(\cdot,b,\tau)}=\hzn{\phi}\cap[-A^{-1}\hzn{\rho}]$.

\paragraph{Part 2} 
For the dual coercivity equivalence,
let $\hat g(u)= g_\tau(u)-\ip{b}{u}$, which is the objective of the reduced
dual~\ref{GeneralValueDualr}. By \eqref{g conj2}, $\hat g^*(0) =
g_\tau^*(b)=\cl{v(\cdot,\tau)}\le v(b,\tau)$.  Therefore, the result
follows from \cite[Corollary 14.2.2]{RTR} since by~\eqref{fdomchar},
$\dom{v(\cdot,\tau)}=\dom{\rho}+A\dom{\phi}$.

\subsection*{Proof of Theorem~\ref{value function sd}}

The expression for $f^*$ is derived in~\eqref{f conj}.
The weak and strong duality relationships as well as the expression for 
$\sd v$ follow immediately
from~\cite[Theorem~11.39]{RTRW}. 

Next, note that
\begin{equation}
\label{fdomchar}
 \dom{f(\cdot,b,\tau)} \neq \emptyset \iff
 \left[
 \begin{array}{c}
 \exists x \in \lev{\phi}{\tau}\\ b - Ax \in \dom{\rho}
 \end{array}
 \right]
\iff
b\in \dom{\rho} + A\lev{\phi}{\tau}.
\end{equation}

Now assume that $b\in\intr{\dom{\rho}+A(\lev{\phi}{\tau})}$.  Recall
from \cite[Theorem 6.6 and Corollary 6.6.2]{RTR} that
 \begin{equation}\label{intr to ri}
 \intr{\dom{\rho}+A(\lev{\phi}{\tau})}=\ri{\dom{\rho}}+A(\ri{\lev{\phi}{\tau}}).
 \end{equation}
 Moreover, by
 \cite[Theorem 7.6]{RTR}, for any convex function $p$,
 \begin{equation}\label{ri of lev}
 \ri{\lev{p}{\tau}}
 =\set{x\in\ri{\dom{p}}}{p(x)< \tau}.
 \end{equation}
 Since $b\in\intr{\dom{\rho}+A(\lev{\phi}{\tau})}$, \eqref{intr to
   ri}--\eqref{ri of lev} imply the existence of
 $\bw\in\ri{\dom{\rho}}$ and $\bx\in\ri{\dom{\phi}}$ with
 $\phi(\bx)<\tau$ such that $b=\bw+A\bx$. Since $\phi$ is relatively
 continuous on the relative interior of its domain \cite[Theorem
 10.1]{RTR}, there exists $\del>0$ such that
 \begin{align*}
   (\bw+\del\uball)\cap\dom{\rho}\subset\ri{\dom{\rho}},
\\ (\bx+\del\uball)\cap\dom{\phi}\subset\ri{\dom{\phi}},
\\ \phi(x)<\half(\phi(\bx)+\tau)\; \forall\,
   x\in (\bx+\del\uball)\cap\dom{\phi}.
 \end{align*}
 Set $S_\rho=(\bw+\del\uball)\cap\dom{\rho}$ and $S_\phi=(\bx+\del\uball)\cap\dom{\phi}$.
 Since
 \begin{align*}
 \cone{S_\rho+AS_\phi-b}&=\cone{S_\rho-\bw}+A\cone{S_\phi-\bx}\\
 &=\Span{\dom{\rho}-\bw}+A\;\Span{\dom{\phi}-\bx}\\
 &=\Span{\dom{\rho}+A\dom{\phi}-b} \\
 &\supset\cone{\dom{\rho}+A\dom{\phi}-b}\\
 &=\Rm\qquad\left(b\in\intr{\dom{\rho}+A(\lev{\phi}{\tau})}\right),
 \end{align*}
 we have $0\in\intr{ S_\rho+AS_\phi-b}$.
 Therefore, there exits an $\eps>0$ such that $b+\eps\uball\subset S_\rho+AS_\phi$.
 Consequently, if $\hat b\in  b+\eps\uball$ and $|\hat\tau-\tau|<\half(\tau-\phi(\bx))$, then
 $\dom{f(\cdot,\hat b,\hat \tau)}\ne\emptyset$ and so $(\hat b,\hat \tau)\in\dom{v}$.

 On the other hand, if $(b,\tau)\in\intr{\dom{v}}$, then $\dom{f(\cdot,\hat b,\hat \tau)}\ne\emptyset$
 for all $(\hat b,\hat \tau)$ near $(b,\tau)$ so that
 $\dom{f(\cdot,\hat b,\tau)}\ne\emptyset$ for all $\hat b$ near $b$. Hence
 $b\in\intr{\dom{\rho}+A(\lev{\phi}{\tau})}$.

\subsection*{Proof of Theorem~\ref{reduced primal-dual}}

\paragraph{Part 1} First note that~\eqref{ReducedOptCond} is
equivalent to the optimality condition
\begin{equation}
\label{zero subdiff}
0 \in -A^T \sd{\rho}(b - A\bx) +  \sd\indicator{\bx}{\lev{\phi}{\tau}}
\end{equation}
for the problem~\ref{ConstrainedClass}, and hence by~\cite[Theorem
23.8]{RTR}, $\bx$ solves~\ref{ConstrainedClass}. Moreover,
by~\cite[Theorem 23.5]{RTR},~\eqref{ReducedOptCond} is equivalent to
\[
b-A\bx \in \sd{\rho^*}{(\bu)}, \quad \bx \in \sd \support{A^T\bu}{\lev{\phi}{\tau}},
\]
or, equivalently,
\begin{equation}
\label{zero subdiff 2}
b \in \sd{\rho^*}{(\bu)} + A\sd\support{A^T\bu}{\lev{\phi}{\tau}},
\end{equation}
which by~\cite[Theorem 23.8]{RTR} implies that $\bu$ solves the reduced dual~\ref{GeneralValueDualr}.

\paragraph{Part 2} If $\bx$ solves~\ref{ConstrainedClass}, then
\[
0 \in \sd{\left[ \rho\big(b - A(\cdot)\big) + \indicator{\cdot}{\lev{\phi}{\tau}}\right]}({\bx}),
\]
which by~\cite[Theorems 23.8, 23.9]{RTR} and~\eqref{primal regularity} is equivalent to~\eqref{zero subdiff},
which in turn is equivalent to~\eqref{ReducedOptCond}.

\paragraph{Part 3} If $\bu$ solves~\ref{GeneralValueDualr}, then
\[
b \in \sd\left[{\rho^*(\cdot) + \support{A^T(\cdot)}{\lev{\phi}{\tau}}}\right](\bu),
\]
which by~\cite[Theorems 23.8, 23.9]{RTR} and~\eqref{dual regularity}
is equivalent to~\eqref{zero subdiff 2}, which in turn is equivalent
to~\eqref{ReducedOptCond}.

\paragraph{Part 4}
The equivalence~\eqref{subdiff explicit 2} follows from~\eqref{subdiff
  explicit}, Part 2 of Lemma~\ref{mu existence}, and the fact that
$A^T\bu\in\ncone{\bx}{\lev{\phi}{\tau}}$ if and only if
$\bx\in\sd\support{A^T\bu}{\lev{\phi}{\tau}}$.

To see~\eqref{subdiff explicit}, note that \eqref{eq:primal
  coercivity}, \eqref{primal regularity}, and Part 1 of Lemma
\ref{coercivity 1} imply that the primal objective is coercive, so a
solution $\bx$ exists.  Hence, by Part 2, there exists $\bu$ so that
$(\bx, \bu)$ satisfies~\eqref{ReducedOptCond}.

Analogously, \eqref{eq:dual coercivity}, \eqref{dual regularity}, and
Part 2 of Lemma \ref{coercivity 1} imply that the solution $\bu$ to
the dual exists, and so by Part 3, there exists $\bx$ such that the
pair $(\bx, \bu)$ satisfies~\eqref{ReducedOptCond}.  In either case,
the subdifferential is nonempty and is given by \eqref{subdiff
  explicit}.

\subsection*{Proof of Lemma~\ref{case 1 details}}

Formula \eqref{gauge=support} is just \cite[Theorem 14.5]{RTR}.
The first equation in \eqref{gauge dom} is obvious and the second follows from
\eqref{gauge=support} and the definition of the barrier cone.
The formula \eqref{c1 dom phi} is now obvious.
Formulas \eqref{c1 lev phi} and \eqref{c1 hzn phi} follow immediately from the
definitions and \cite[Corollary 8.3.3]{RTR}.
Formula \eqref{c1 bar lev phi} follows from \eqref{c1 hzn phi}, \cite[Corollary 14.2.1]{RTR},
and \cite[Corollary 16.4.2]{RTR}.

First note that \eqref{c1 cq} implies that
$\ri{\cone{U}}\cap\ri{X}\ne\emptyset$.
Hence, the formula \eqref{c1 conj} follows from
\cite[Theorem 16.4]{RTR} and \eqref{c1 dom phi}.
To see \eqref{c1 conj pi}, observe that the expression on the RHS is again an
infimal convolution for which $\inf=\min$ for the same reason as for \eqref{c1 conj}.
The equivalence with $(\phi^*)^\pi(z,\mu)$ follows from the calculus rules in
section~\ref{sec:calculus}.
For formula \eqref{c1 support lev phi}, first note that
\begin{eqnarray*}
{\inf_{\mu\ge 0}[\tau\mu +(\phi^*)^\pi(z,\mu)]}&=&
\inf_{\mu\ge 0}\left[\tau\mu +\inf_s[\support{z-s}{X}+\indicator{s}{\mu U^\circ}]\right]\\
&=&
\inf_s\left[\support{z-s}{X}+\inf_{\mu\ge 0}[\tau\mu +\indicator{s}{\mu U^\circ}]\right]\\
&=&
\inf_s\left[\support{z-s}{X}+\tau\gauge{s}{U^\circ}\right].
\end{eqnarray*}
Again, the final infimum in this derivation is an infimal convolution for which
$\inf=\min$ for the same reasons as in \eqref{c1 conj} since, by \eqref{c1 dom phi}
and \cite[Theorem 14.5]{RTR},
\[
\dom{\left((\tau\gauge{\cdot}{U^\circ})^*\right)}=\dom{\left((\support{\cdot}{\tau U})^*\right)}
=\dom{\indicator{\cdot}{\tau U}}=\tau U.
\]
Therefore, an optimal $\bs$ in this infimal convolution exists giving
$\bmu=\gauge{\bs}{U^\circ}$ as the optimal solution to
the first $\min$ in \eqref{c1 support lev phi}.

Formula \eqref{c1 ncone lev phi} is an immediate consequence of
\eqref{c1 lev phi}, \eqref{c1 cq}, and \cite[Corollary 23.8.1]{RTR}.

\subsection*{Proof of Theorem~\ref{case 1 primal-dual}}

By~\eqref{c1 lev phi} and the calculus rules for the relative
interior~\cite[Section 6]{RTR}, \eqref{primal regularity} and
\eqref{c1 primal regularity} are equivalent.  Similary, by~\eqref{c1
  bar lev phi} and \cite[Theorem 6.3]{RTR}, \eqref{dual regularity}
and~\eqref{c1 dual regularity} are equivalent.

\paragraph{Part 1} Since \eqref{c1 cq} holds, the formula \eqref{c1 ncone lev phi} holds and so
\eqref{cases1ab} and \eqref{ReducedOptCond} are equivalent. Hence, the
result follows from Part 1~of Theorem \ref{reduced primal-dual}.

\paragraph{Part 2} Since \eqref{primal regularity} and \eqref{c1 primal regularity} are equivalent,
the result follows from Part 2~of Theorem \ref{reduced primal-dual}.

\paragraph{Part 3} Since \eqref{dual regularity} and \eqref{c1 dual regularity} are equivalent,
the result follows from Part 3~of Theorem \ref{reduced primal-dual}.

\paragraph{Part 4}
By \eqref{c1 hzn phi}, \eqref{c1 coercive primal} is equivalent
to~\eqref{eq:primal coercivity} and~\eqref{primal regularity}, and,
by~\eqref{c1 dom phi}, \eqref{c1 coercive dual} is equivalent
to~\eqref{eq:dual coercivity} and~\eqref{dual regularity}.  Therefore,
by Theorem \ref{reduced primal-dual}, \eqref{c1 subdiff explicit} is
equivalent to \eqref{subdiff explicit} since
$\tau\gauge{s}{U^\circ}=\inf_{\mu\ge 0}[\tau\mu+\indicator{s}{\mu
  U^\circ}$.  The final equivalence is identical to that of Theorem
\ref{reduced primal-dual}.

\subsection*{Proof of Lemma~\ref{c2 phi}}

The formula for $\dom\phi$ follows from~\eqref{LQ gauge}.
Indeed, by~\eqref{LQ gauge}, $x\in\dom{\phi}$ if and only if there exists
$s\in\R^k$ such that $x-Ls\in\dom{\gauge{\cdot}{U^\circ}}=\cone{U^\circ}$, 
or equivalently, $x\in \cone{U^\circ}+\Ran{L}=\cone{U^\circ}+\Ran{B}$.
The formula for $\hzn{\phi}$ follows immediately from \cite[Theorem 14.2]{RTR}
and \eqref{case 2 phi conj}.
In particular, $\phi$ is coercive if and only if $\{0\}=\hzn{\phi}$, or equivalently,
$\cone{U}=\R^n$, i.e., $0\in\intr{U}$.

Next we show that
the $\lam$ given
  in \eqref{minimizingLambda} solves  \eqref{ind conj inf rep}. First observe
  that the optimal $\lambda$ must be greater than $\gauge{w}{U}$, and
  from elementary calculus, the minimizer of the hyperbola
  $\frac{1}{2\lambda}\|w\|_B^2 + \tau\lambda $ for $\lam\ge 0$ is
  given by $\|w\|_B/\sqrt{2\tau}$.  Therefore, the minimizing $\lambda$
  is given by \eqref{minimizingLambda}.  Substituting this value into
  \eqref{ind conj inf rep} gives~\eqref{indicatorConjugate}.

  It is now easily shown that the function in~\eqref{indicatorConjugate} is
  lower semi-continuous. Therefore, the equivalence in \eqref{ind conj inf rep}
  follows from \eqref{level support}.

\subsection*{Proof of Theorem~\ref{case 2 sd v}}

By \cite[Theorem 7.6]{RTR},
\[\ri{\lev{\phi}{\tau}}=\set{x}{x\in\ri{\dom{\phi}},\ \phi(x)<\tau}.\]
Hence, by Lemma \ref{c2 phi}, the equivalence between
\eqref{eq:reduced primal-dual}
and \eqref{c2 primal regularity}, \eqref{c2 dual regularity},
\eqref{c2 primal existence}, \eqref{c2 dual existence}, respectively,
is easily seen.  Therefore, Parts 1--4 follow immediately from
Theorem~\ref{reduced primal-dual}.

\subsection*{Proof of Corollary~\ref{c2 multiplier cor}}

Condition~\eqref{c2 multiplier rule a} occurs when $\bmu=0$ since
$0^+\sd\phi(\bx)=\ncone{\bx}{\dom\phi}$. When $\bmu>0$,
by~\cite[Theorem 23.5]{RTR}, $\partial \phi(x)=\argmax_{w\in
  U}[\ip{x}{w}-\half \ip{w}{Bw}]$, so that $w\in\partial \phi(x)$ if
and only if $x\in Bw+\ncone{w}{U}$.

\section*{Acknowledgments} The authors are grateful to two anonymous
referees for their remarkably careful reading of an intricate
paper. Their detailed list of comments and suggestions led to many
fixes and improvements.

\bibliographystyle{plain}
\bibliography{biblio}

\begin{thebibliography}{10}

\bibitem{Anderson:1979}
B.~D.~O. Anderson and J.~B. Moore.
\newblock {\em Optimal Filtering}.
\newblock Prentice-Hall, Englewood Cliffs, N.J., USA, 1979.

\bibitem{Aravkin2011}
A.~Aravkin, B.~Bell, J.~V. Burke, and G.~Pillonetto.
\newblock An $\ell_1$-{L}aplace robust {K}alman smoother.
\newblock {\em IEEE Transactions on Automatic Control}, 2011.

\bibitem{SYSID2012tks}
A.~Aravkin, J.~Burke, and G.~Pillonetto.
\newblock Robust and trend following kalman smoothers using {S}tudent's t.
\newblock In {\em Proc. of SYSID}, 2012.

\bibitem{ABP:2013}
A.~Aravkin, James Burke, and Gianluigi Pillonetto.
\newblock Sparse/robust estimation and kalman smoothing with nonsmooth
  log-concave densities: Modeling, computation, and theory.
\newblock {\em To appear in the Journal of Machine Learning Research}, Accepted
  April 2013.

\bibitem{AravkinFHV:2012}
A.~Aravkin, M.~P. Friedlander, F.~Herrmann, and T.~van Leeuwen.
\newblock Robust inversion, dimensionality reduction, and randomized sampling.
\newblock {\em Mathematical Programming}, 134(1):101--125, 2012.

\bibitem{Aronszajn}
N.~Aronszajn.
\newblock Theory of reproducing kernels.
\newblock {\em Trans. of the American Mathematical Society}, 68:337--404, 1950.

\bibitem{AusTeb:2003}
A.~Auslender and M.~Teboulle.
\newblock {\em Asymptotic Cones and Functions in Optimization and Variational
  Inequalities}.
\newblock Springer, 2003.

\bibitem{BauCom:2011}
H.~H. Bauschke and P.~L. Combettes.
\newblock {\em Convex Analysis and Monotone Operator Theory in Hilbert Spaces}.
\newblock Springer, 2011.

\bibitem{BergFrie:2007b}
E.~{van den} Berg and M.~P. Friedlander.
\newblock {SPGL1}: A solver for large-scale sparse reconstruction.
\newblock Available at
  \url{http://www.cs.ubc.ca/labs/scl/index.php/Main/Spgl1}, June 2007.

\bibitem{BergFriedlander:2008}
E.~{van den} Berg and M.~P. Friedlander.
\newblock Probing the pareto frontier for basis pursuit solutions.
\newblock {\em SIAM Journal on Scientific Computing}, 31(2):890--912, 2008.

\bibitem{BergFrie:2010}
E.~{van den} Berg and M.~P. Friedlander.
\newblock Theoretical and empirical results for recovery from multiple
  measurements.
\newblock {\em IEEE Transactions on Information Theory}, 56(5):2516--2527,
  2010.

\bibitem{BergFriedlander:2011}
E.~{van den} Berg and M.~P. Friedlander.
\newblock Sparse optimization with least-squares constraints.
\newblock {\em SIAM J. Optimization}, 21(4):1201--1229, 2011.

\bibitem{Borwein00}
J.~M. Borwein and A.~S. Lewis.
\newblock {\em Convex Analysis and Nonlinear Optimization}.
\newblock Springer, 2000.

\bibitem{Boyd04}
S.~Boyd and L.~Vandenberghe.
\newblock {\em Convex Optimization}.
\newblock Cambridge University Press, 2004.

\bibitem{Brockett}
R.~Brockett.
\newblock {\em Finite Dimensional Linear Systems}.
\newblock John Wiley and Sons, Inc., 1970.

\bibitem{ComLagThi:1994}
C.~Combari, M.~Laghdir, and L.~Thibault.
\newblock Sous-diff$\acute{\mathrm{e}}$rentiels de fonctions convexes
  compos$\acute{\mathrm{e}}$es.
\newblock {\em Ann. Sci. Math. Qu$\acute{\mathrm{e}}$bec}, 18(2):119--148,
  1994.

\bibitem{Cucker}
F.~Cucker and S.~Smale.
\newblock On the mathematical foundations of learning.
\newblock {\em Bulletin of the American mathematical society}, 39:1--49, 2001.

\bibitem{DON2005Ta}
D.~L. Dohono and J.~Tanner.
\newblock Sparse nonnegative solution of underdetermined linear equations by
  linear programming.
\newblock {\em Proc. Nat. Acad. Sci.}, 102(27):9446--9451, 2005.

\bibitem{Donoho2006}
D.~Donoho.
\newblock Compressed sensing.
\newblock {\em IEEE Trans. on Information Theory}, 52(4):1289--1306, 2006.

\bibitem{LARS2004}
B.~Efron, T.~Hastie, L.~Johnstone, and R.~Tibshirani.
\newblock Least angle regression.
\newblock {\em Annals of Statistics}, 32:407--499, 2004.

\bibitem{EkeTem:1976}
I.~Ekeland and R.~Temam.
\newblock {\em Convex Analysis and Variational Problems}.
\newblock Elsevier, 1976.

\bibitem{Evgeniou99}
T.~Evgeniou, M.~Pontil, and T.~Poggio.
\newblock Regularization networks and support vector machines.
\newblock {\em Advances in Computational Mathematics}, 13:1--150, 2000.

\bibitem{Farahmand2011}
S.~Farahmand, G.B. Giannakis, and D.~Angelosante.
\newblock Doubly robust smoothing of dynamical processes via outlier sparsity
  constraints.
\newblock {\em IEEE Transactions on Signal Processing}, 59:4529--4543, 2011.

\bibitem{Gao2008}
J.~Gao.
\newblock Robust l1 principal component analysis and its {B}ayesian variational
  inference.
\newblock {\em Neural Computation}, 20(2):555--572, February 2008.

\bibitem{Hastie90}
T.~J. Hastie and R.~J. Tibshirani.
\newblock Generalized additive models.
\newblock In {\em Monographs on Statistics and Applied Probability}, volume~43.
  Chapman and Hall, London, UK, 1990.

\bibitem{Hastie01}
T.~J. Hastie, R.~J. Tibshirani, and J.~Friedman.
\newblock {\em The Elements of Statistical Learning. Data Mining, Inference and
  Prediction}.
\newblock Springer, Canada, 2001.

\bibitem{HFY:2012}
F.~J. Herrmann, M.~P. Friedlander, and O.~Yilmaz.
\newblock Fighting the curse of dimensionality: compressive sensing in
  exploration seismology.
\newblock {\em IEEE Signal Proc. Magazine}, 29(3):88--100, 2012.

\bibitem{Hiriart-Urruty01}
J.~B. Hiriart-Urruty and C.~Lemar\'echal.
\newblock {\em Fundamentals of Convex Analysis}.
\newblock Springer, New York, NY, USA, 2001.

\bibitem{Huber81}
P.~J. Huber.
\newblock {\em Robust Statistics}.
\newblock Wiley, New York, NY, USA, 1981.

\bibitem{MacKay}
D.~J.~C. MacKay.
\newblock {B}ayesian interpolation.
\newblock {\em Neural Computation}, 4:415--447, 1992.

\bibitem{McKayARD}
{D. J. C.} Mackay.
\newblock {B}ayesian non-linear modelling for the prediction competition.
\newblock {\em ASHRAE Trans.}, 100(2):3704--3716, 1994.

\bibitem{Mar:1987}
H.~M. Markowitz.
\newblock {\em Mean-Variance Analysis in Portfolio Choice and Capital Markets}.
\newblock Frank J.~Fabozzi Associates, New Hope, Penn., 1987.

\bibitem{Pontil98}
M.~Pontil and A.~Verri.
\newblock Properties of support vector machines.
\newblock {\em Neural Computation}, 10:955--974, 1998.

\bibitem{RTR}
R.~T. Rockafellar.
\newblock {\em Convex Analysis}.
\newblock Princeton Landmarks in Mathematics. Princeton University Press, 1970.

\bibitem{Rockafellar:1993}
R.~T. Rockafellar.
\newblock Lagrange multipliers and optimality.
\newblock {\em SIAM Review}, 35(2):pp. 183--238, 1993.

\bibitem{RTRW}
R.~T. Rockafellar and R.~J.-B. Wets.
\newblock {\em Variational Analysis}, volume 317.
\newblock Springer, 1998.

\bibitem{Roweis99}
S.~Roweis and Z.~Ghahramani.
\newblock A unifying review of linear {G}aussian models.
\newblock {\em Neural Computation}, 11:305--345, 1999.

\bibitem{Saitoh}
S.~Saitoh.
\newblock {\em Theory of reproducing kernels and its applications}.
\newblock Longman, 1988.

\bibitem{Scholkopf00}
B.~Sch\"{o}lkopf, A.~J. Smola, R.~C. Williamson, and P.~L. Bartlett.
\newblock New support vector algorithms.
\newblock {\em Neural Computation}, 12:1207--1245, 2000.

\bibitem{Tap:2013}
R.~Tapia.
\newblock The isoperimetric problem revisited: extracting a short proof of
  sufficiency from {E}uler's approach to necessity.
\newblock Technical report, Rice University, CAAM, 2013.

\bibitem{Tib96}
R.~Tibshirani.
\newblock Regression shrinkage and selection via the {L}asso.
\newblock {\em J. R. Statist. Soc. B.}, 58(1):267--288, 1996.

\bibitem{Tipping2001}
M.~Tipping.
\newblock Sparse {B}ayesian learning and the relevance vector machine.
\newblock {\em Journal of Machine Learning Research}, 1:211--244, 2001.

\bibitem{Vapnik98}
V.~Vapnik.
\newblock {\em Statistical Learning Theory}.
\newblock Wiley, New York, NY, USA, 1998.

\bibitem{Wipf_IEEE_TIT_2011}
{D. P.} Wipf, {B.D.} Rao, and S.~Nagarajan.
\newblock Latent variable {B}ayesian models for promoting sparsity.
\newblock {\em IEEE Transactions on Information Theory (to appear)}, 2011.

\bibitem{Zal:2002}
C.~Z$\breve{\mathrm{a}}$linescu.
\newblock {\em Convex Analysis in General Vector Spaces}.
\newblock World Scientific, 2002.

\end{thebibliography}

\end{document}